\documentclass[12pt,a4paper]{article}%
\usepackage{amsmath}
\usepackage{graphicx}
\usepackage{a4}%
\usepackage{amsfonts}%
\usepackage{amssymb}

\begin{document}

\title{Multiplicative rule of Schubert classes}
\author{Haibao Duan\\Institute of Mathematics, Chinese Academy of Sciences\\Beijing 100080, dhb@math.ac.cn}
\date{\ \ \ }
\maketitle

\begin{abstract}
Let $G$ be a compact connected Lie group and $H$, the centralizer of a
one-parameter subgroup in $G$. Combining the ideas of Bott-Samelson
resolutions of Schubert varieties and the enumerative formula on a twisted
product of $2$ spheres obtained in [Du$_{2}$], we obtain an explicit formula
for multiplying Schubert classes in the flag manifold $G\diagup H$.

\begin{description}
\item 2000 Mathematical Subject Classification: 14N15 (14M10).

\item Key words and phrases: Schubert varieties, interesection multiplicities,
Cartan numbers

\end{description}
\end{abstract}

\section{\textbf{Introduction}\ \ \ }

Let $G$ be a compact connected Lie group and $H$, the centralizer of a
one-parameter subgroup in $G$. The Weyl of $G$ (resp. of $H$) is denoted by
$W$ (resp. $W^{\prime}$). The set $W/W^{\prime}$ of left cosets of $W^{\prime
}$ in $W$ can be identified with the subset of $W$:

$\qquad\qquad\overline{W}=\{w\in W\mid l(w_{1})\geq l(w)$ for all $w_{1}\in
wW^{\prime}\}$,\

\noindent where $l:W\rightarrow\mathbb{Z}$ is the length function relative to
a fixed maximal torus $T$ in $G$.

It is known from Bruhat-Chevalley that the flag manifold $G/H=\{gH\mid g\in
G\}$ admits a canonical decomposition into cells, indexed by elements of
$\overline{W}$,

$\qquad G/H=\underset{w\in\overline{W}}{\cup}X_{w}(H)$, $\quad\dim
X_{w}=2l(w)$,

\noindent with each cell $X_{w}(H)$ the closure of an algebraic affine space,
known as a \textsl{Schubert variety} in $G/H$ [BGG]. Since only even
dimensional cells are involved, the set of fundamental classes\textsl{
}$[X_{w}(H)]\in H_{2l(w)}(G/H)$, $w\in\overline{W}$,\textsl{ }form an additive
basis of the homology $H_{\ast}(G/H)$. The cocycle class $P_{w}(H)\in
H^{2l(w)}(G/H)$, $w\in\overline{W}$, defined by the Kronecker pairing as
$\left\langle P_{w}(H),[X_{u}(H)]\right\rangle =\delta_{w,u}$, $w,u\in
\overline{W}$, is called the \textsl{Schubert class corresponding to }$w$.
Clearly one has

\textbf{Basis Theorem.} \textsl{The set of Schubert classes }$\{P_{w}(H)\mid
$\textsl{\ }$w\in\overline{W}\}$\textsl{\ constitutes an additive basis for
the cohomology }$H^{\ast}(G/H)$\textsl{.}

One immediate consequence is that the product of two arbitrary Schubert
classes can be expressed in terms of Schubert classes. Precisely, given
$u,v\in\overline{W}$, one has the expression

$\qquad P_{u}(H)\cdot P_{v}(H)=\sum\limits_{l(w)=l(u)+l(v),w\in\overline{W}%
}a_{u,v}^{w}P_{w}(H)$, $a_{u,v}^{w}\in\mathbb{Z}$

\noindent in $H^{\ast}(G/H)$. Since the Chow ring $A^{\ast}(G/H)$ is
canonically isomorphic to the integral cohomology $H^{\ast}(G/H)$, the
following Problem is of fundamental importance in the intersection theory of
$G/H$.

\textbf{Problem.} \textsl{Find the number }$a_{u,v}^{w}$ \textsl{for given
}$w,u,v\in\overline{W}$\textsl{,} $l(w)=l(u)+l(v)$\textsl{.}

\bigskip

If $G$ is the unitary group $U(n)$ of rank $n$ and $H=U(k)\times U(n-k)$, the
flag manifold $G/H$ is the Grassmannian $G_{n,k}$ of $k$-planes through the
origin in $\mathbb{C}^{n}$. In this case, a combinatorial description for
$a_{u,v}^{w}$ is given by the \textsl{Littlewood-Richardson rule},\textsl{
}one cornerstone of the Schubert calculus for $G_{n,k}$ [S]. It was first
stated by Littlewood and Richardson\textsl{ }in 1934 [LR]. Complete proofs
appeared only in the 1970s (see ``Note and references'' in [M, p.148]).

Another special case is when $H=T$ (a maximal torus in $G$) and if either
$l(u)=1$ or $l(v)=1$. The number $a_{u,v}^{w}$ is seen as certain Cartan
number of $G$ from the \textsl{Chevalley formula}. Chevalley announced the
formula at the end of his address at the 1958 ICM in Edinburgh [Ch$_{1}$],
while a proof was given in his famous manuscript [Ch$_{2}$]. Although
[Ch$_{2}$] remained unpublished until 1994, this formula became part of the
official literature after the publications of [BGG] by Bernstein et al in
1973, and [De$_{2}$] by Demazure in 1974, where both authors verified it using
different methods (cf. introduction to [Ch$_{2}$] by Borel).

In recent years, inspired by theory of Schubert polynomials of Lascoux and
Sch\"{u}tzenberger [LS], many achievements have been made in generalizing the
classical \textsl{Pieri formula}, which handles the problem for the case where
$G$ is a matrix group and where one of $P_{u}$ and $P_{v}$ is a
\textsl{special Schubert class}. (See [FP, Section 9.10] for more recent
progresses and relevant references).

While the problem in its natural generality remains unsolved\footnote{We quote
from Fulton and Pragacz [FP]: \textsl{there is no \ analogue of the Littlewood
Richardson rule for explicitly multiplying Schubert classes in a flag
manifold}; from Sottile [S]: \textsl{the analog of the Littlewood Richardson
rule is not known for most other flag variety G/P}.
\par
For the cases of matrix groups, the theory of Schubert polynomials was
developed to make explicit computation with Schubert classes possible (cf.
introduction to [BH]).}, the further problem of determining the multiplicative
rule of Schubert classes in the quantum cohomology of $G/H$ has appeared on
the agenda, where the analogue of the coefficients $a_{u,v}^{w}$ are known as
\textsl{Gromov-Witten} \textsl{numbers} (cf. [FP, p.134], [CF]).

It was announced in [Du$_{2}$] that, combining the ideas of Bott-Samelson
resolutions of Schubert varieties and the enumerative formula on a twisted
products of $2$-spheres obtained in [Du$_{2}$], it is possible to find a
unified formula that expresses $a_{u,v}^{w}$ in terms of certain Cartan
numbers of $G$. This paper is devoted to complete this project.

\section{\bigskip Main result}

A few notations will be needed in presenting our result. Throughout this paper
$G$ is a compact connected Lie group with a fixed maximal torus $T$. We set
$n=\dim T$.

\smallskip\ \ \

\textbf{2.1. Geometry of Cartan subalgebra. }Equip the Lie algebra $L(G)$ of
$G$ with an inner product $(,)$ so that the adjoint representation acts as
isometries of $L(G)$. The\textsl{\ Cartan subalgebra} of $G$ is the Euclidean
subspace $L(T)$ of $L(G)$.

The restriction of the exponential map $\exp:L(G)\rightarrow G$ to $L(T)$
defines a set $D(G)$ of $m=\frac{1}{2}(\dim G-n)$ hyperplanes in $L(T)$, i.e.
the set of\textsl{\ singular hyperplanes }through the origin in $L(T)$. These
planes divide $L(T)$ into finitely many convex cones, called the \textsl{Weyl
chambers} of $G$. The reflections $\sigma$ of $L(T)$ in the these planes
generate the Weyl group $W$ of $G$. Let $\Phi\subset L(T)$ be the \textsl{root
system associated to $W$ }([Hu]). Recall that if $\beta,\beta^{\prime}\in\Phi
$, \textsl{the Cartan number} $\beta\circ\beta^{\prime}=2(\beta,\beta^{\prime
})/(\beta^{\prime},\beta^{\prime})$ is an integer (only $0,\pm1,\pm2,\pm3$ can occur).\

Fix, once and for all, a regular point $\alpha\in L(T)\backslash\underset{L\in
D(G)}{\cup}L$, and let $\Phi^{+}$ (resp. $\Delta$) be the set of positive
roots (resp. simple roots) relative to $\alpha$. The set $D(G)$ can now be
indexed by $\Phi^{+}$ as $\{L_{\beta}\mid\beta\in\Phi^{+}\}$, where $L_{\beta
}$ is the singular plane corresponding to the root $\beta$. For a $\beta
\in\Phi^{+}$, write $\sigma_{\beta}\in W$ for the reflection of $L(T)$ in
$L_{\beta}$. If $\beta\in\Delta$ we call $\sigma_{\beta}$ a \textsl{simple
reflection}.

\bigskip

It is known that the set of simple reflections $\{\sigma_{\beta}\mid\beta
\in\Delta\}$ generates $W$. That is, any $w\in W$ admits a factorization of
the form

\noindent(2.1) $\qquad\qquad w=\sigma_{\beta_{1}}\circ\cdots\circ\sigma
_{\beta_{k}}$, $\beta_{i}\in\Delta$.

\textbf{Definition 1}. The \textsl{length} $l(w)$ of an $w\in W$ is the least
number of factors in all decompositions of $w$ in the form (2.1). The
decomposition (2.1) is said \textsl{reduced} if $k=l(w)$.

If (2.1) is a reduced decomposition, the $k\times k$ (strictly upper
triangular) matrix $A_{w}=(a_{i,j})$ with

$\qquad\qquad a_{i,j}=\{%
\begin{array}
[c]{c}%
0\text{ if }i\geq j\text{;\qquad}\\
-\beta_{i}\circ\beta_{j}\text{ if }i<j
\end{array}
$

\noindent is called \textsl{the Cartan matrix of $w$ associated to the
decomposition }(2.1).

\bigskip

\textbf{Example 1. }By resorting to the geometry of the Cartan subalgebra
$L(T)$ there is geometric method to find a reduced decomposition (hence a
Cartan matrix) of an $w\in W$.

Picture $W$ as the $W$-orbit $\{w(\alpha)\in L(T)\mid w\in W\}$ of the regular
point $\alpha$. Given an $w\in W$ let $C_{w}$ be a straight line segment in
$L(T)$ from the Weyl chamber containing $\alpha$ to $w(\alpha)$ that crosses
the planes in $D(G)$ one at a time. Assume that they are met in the order
$L_{\alpha_{1}},\cdots,L_{\alpha_{k}}$, $\alpha_{i}\in\Phi^{+}$. We have
$l(w)=k$ and $w=\sigma_{\alpha_{k}}\circ\cdots\circ\sigma_{\alpha_{1}}$ (cf.
[Du$_{2}$] or [Han]). Set

$\qquad\beta_{1}=\alpha_{1},\quad\beta_{2}=\sigma_{\alpha_{1}}(\alpha
_{2}),\cdots,\quad\beta_{k}=\sigma_{\alpha_{1}}\circ\cdots\circ\sigma
_{\alpha_{k-1}}(\alpha_{k})$.

\noindent Then $\beta_{i}\in\Delta$. Moreover, from

$\qquad\sigma_{\beta_{i}}=\sigma_{\alpha_{1}}\circ\cdots\circ\sigma
_{\alpha_{i-1}}\circ\sigma_{\alpha_{i}}\circ\sigma_{\alpha_{i-1}}\circ
\cdots\circ\sigma_{\alpha_{1}}$,

\noindent one verifies easily that $w=\sigma_{\beta_{1}}\circ\cdots\circ
\sigma_{\beta_{k}}$. This decomposition is reduced because of $\beta_{i}%
\in\Delta$ and $l(w)=k$.

As an example of the method consider the case of $G_{2}$, the exceptional Lie
group of rank $2$. With $\dim L(T)=2$, the singular lines, denoted by $L_{i}$,
$i\leq6$, are depicted in the figure. Taking a regular point $\alpha\in L(T)$
as marked, the set of simply roots is $\Delta=\{\beta_{1},\beta_{2}\}$. Let
$W$ be the Weyl group of $G_{2}$.%

\begin{figure}
[h]
\begin{center}
\includegraphics[
trim=0.000000in 0.000000in -0.049297in -0.352822in,
height=3.4982in,
width=4.2921in
]%
{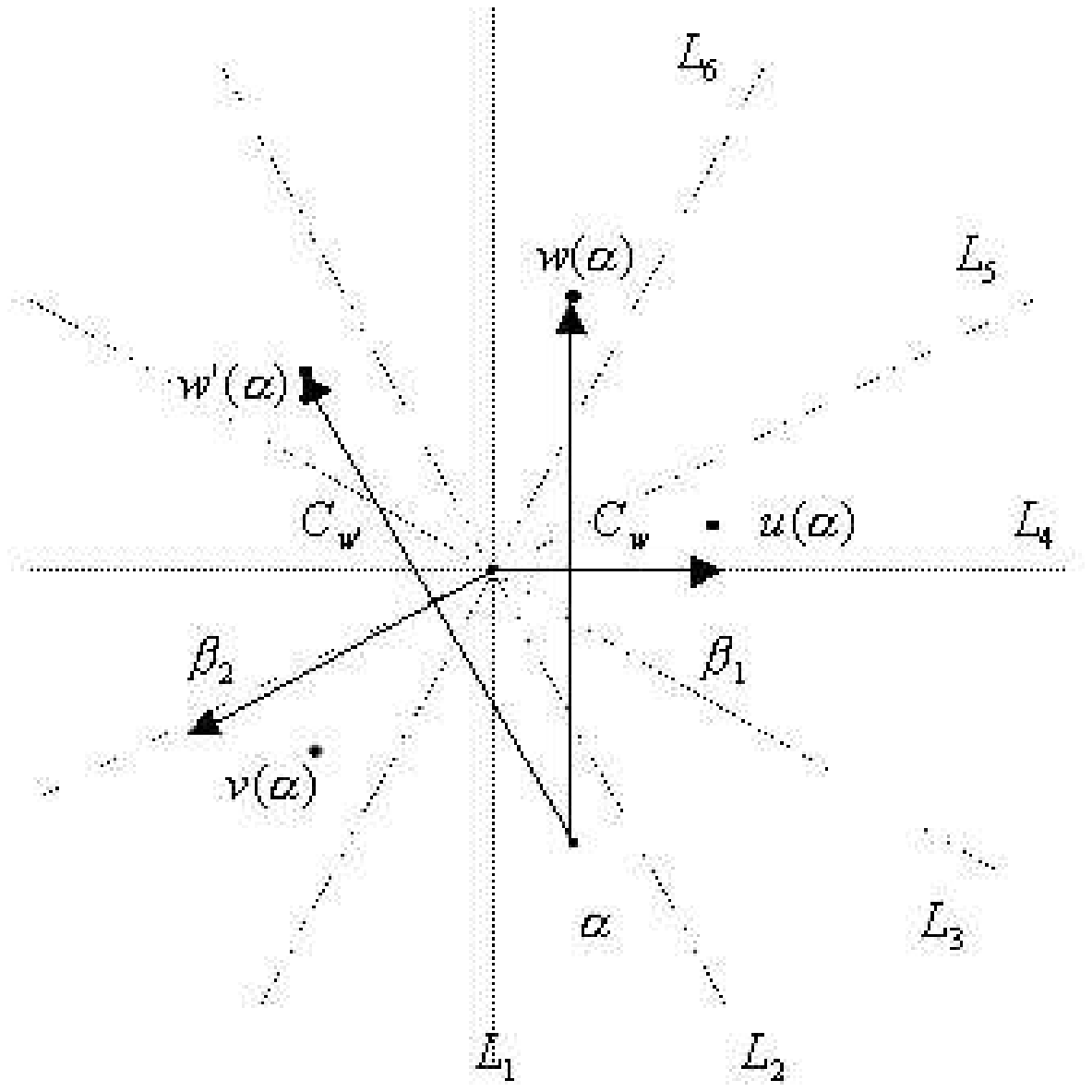}%
\end{center}
\end{figure}

For the elements $w,w\prime\in W$ specified by the vectors $w(\alpha
),w\prime(\alpha)\in L(T)$ in the figure, we get from the segments $C_{w}$ and
$C_{w^{\prime}}$ the reduced decompositions

$\qquad w=\sigma_{\beta_{2}}\circ\sigma_{\beta_{1}}\circ\sigma_{\beta_{2}%
}\circ\sigma_{\beta_{1}}\circ\sigma_{\beta_{2}}$;

$\qquad w\prime=\sigma_{\beta_{1}}\circ\sigma_{\beta_{2}}\circ\sigma
_{\beta_{1}}\circ\sigma_{\beta_{2}}\circ\sigma_{\beta_{1}}$.

\noindent From the Cartan matrix of $G_{2}$ ([Hu, p.59])

$\qquad\left(
\begin{array}
[c]{cc}%
\beta_{1}\circ\beta_{1} & \beta_{1}\circ\beta_{2}\\
\beta_{2}\circ\beta_{1} & \beta_{2}\circ\beta_{2}%
\end{array}
\right)  =\left(
\begin{array}
[c]{cc}%
2 & -1\\
-3 & 2
\end{array}
\right)  $

\noindent we read that the Cartan matrices of $w$ and $w^{\prime}$ associated
to the decompositions are respectively

$A_{w}=\left(
\begin{array}
[c]{ccccc}%
0 & 3 & -2 & 3 & -2\\
0 & 0 & 1 & -2 & 1\\
0 & 0 & 0 & 3 & -2\\
0 & 0 & 0 & 0 & 1\\
0 & 0 & 0 & 0 & 0
\end{array}
\right)  $ and $A_{w^{\prime}}=\left(
\begin{array}
[c]{ccccc}%
0 & 1 & -2 & 1 & -2\\
0 & 0 & 3 & -2 & 3\\
0 & 0 & 0 & 1 & -2\\
0 & 0 & 0 & 0 & 3\\
0 & 0 & 0 & 0 & 0
\end{array}
\right)  $.

\bigskip\ \ \ \

\textbf{2.2. The triangular operators.} Let $\mathbb{Z}[x_{1},\cdots
,x_{k}]=\oplus_{r\geq0}\mathbb{Z}[x_{1},\cdots,x_{k}]^{(r)}$ be the ring of
integral polynomials in $x_{1},\cdots,x_{k}$, graded by $\mid x_{i}\mid=1$.

\textbf{Definition 2. }Given an $k\times k$ strictly upper triangular integer
matrix $A=(a_{i,j})$ the \textsl{triangular operator} associated to $A$ is the
homomorphism $T_{A}:$ $\mathbb{Z}[x_{1},\cdots,x_{k}]^{(k)}\rightarrow
\mathbb{Z}$ defined recursively by the following \textsl{elimination laws}.

1) if $h\in\mathbb{Z}[x_{1},\cdot\cdot\cdot,x_{k-1}]^{(k)}$, then $T_{A}(h)=0$;

2) if $k=1$ (consequently $A=(0)$), then $T_{A}(x_{1})=1$;

3) if $h\in\mathbb{Z}[x_{1},\cdot\cdot\cdot,x_{k-1}]^{(k-r)}$ with $r\geq1$, then

$\qquad\qquad T_{A}(hx_{k}^{r})=T_{A^{\prime}}(h(a_{1,k}x_{1}+\cdots
+a_{k-1,k}x_{k-1})^{r-1})$,

\noindent where $A^{\prime}$ is the ($(k-1)\times(k-1)$ strictly upper
triangular) matrix obtained from $A$ by deleting the $k^{th}$ column and the
$k^{th}$ row.

By additivity, $T_{A}$ is defined for every $f\in\mathbb{Z}[x_{1},\cdots
,x_{k}]^{(k)}$ using the unique expansion $f=\Sigma h_{r}x_{k}^{r}$ with
$h_{r}\in\mathbb{Z}[x_{1},\cdot\cdot\cdot,x_{k-1}]^{(k-r)}$.

\bigskip

\textbf{Example 2.} Definition 2 gives an effective algorithm to evaluate
$T_{A}$.

For $k=2$ and $A_{1}=\left(
\begin{array}
[c]{cc}%
0 & a\\
0 & 0
\end{array}
\right)  $, then $T_{A_{1}}:$ $\mathbb{Z}[x_{1},x_{2}]^{(2)}\rightarrow
\mathbb{Z}$ is given by

$\qquad\qquad T_{A_{1}}(x_{1}^{2})=0$,

$\qquad\qquad T_{A_{1}}(x_{1}x_{2})=T_{A_{1}^{\prime}}(x_{1})=1$ and

$\qquad\qquad T_{A_{1}}(x_{2}^{2})=T_{A_{1}^{\prime}}(ax_{1})=a$.

\bigskip

For $k=3$ and $A_{2}=\left(
\begin{array}
[c]{ccc}%
0 & a & b\\
0 & 0 & c\\
0 & 0 & 0
\end{array}
\right)  $, then $A_{2}^{\prime}=A_{1}$ and $T_{A_{2}}:$ $\mathbb{Z}%
[x_{1},x_{2},x_{3}]^{(3)}\rightarrow\mathbb{Z}$ is given by

$\qquad T_{A_{2}}(x_{1}^{r_{1}}x_{2}^{r_{2}}x_{3}^{r_{3}})=\{%
\begin{array}
[c]{c}%
0\text{, if }r_{3}=0\text{ and\qquad\qquad\qquad\qquad\qquad}\\
T_{A_{1}}(x_{1}^{r_{1}}x_{2}^{r_{2}}(bx_{1}+cx_{2})^{r_{3}-1}),\text{ if
}r_{3}\geq1\text{,\ }%
\end{array}
$

\noindent where $r_{1}+r_{2}+r_{3}=3$, and where $T_{A_{1}}$ is calculated in
the above.

\bigskip

It is straightforward from Definition 2 that

\textbf{Corollary 1.} \textsl{We have }$T_{A}(x_{1}\cdots x_{k})=1$\textsl{ and}

$\qquad\qquad T_{A}(x_{1}^{r_{1}}\cdots x_{k}^{r_{k}})=0$

\noindent\textsl{whenever }$r_{1}+\cdots+r_{i}>i$\textsl{ for some }$1\leq
i<k$\textsl{.}

\bigskip

The operator $T_{A}$ can also be given by explicit formula. For a sequence
$(r_{1},\cdots,r_{k})$ of $k$ non-negative integers with $\Sigma r_{i}=k$, let
$C(r_{1},\cdots,r_{k})$ be the set of all strictly upper triangular $k\times
k$ matrices $C=(c_{i,j})$ with non-negative integer entries satisfying

$\qquad\underset{s}{\Sigma}c_{s,i}=r_{i}-1+\underset{j}{\Sigma}c_{i,j}$,
$1\leq i\leq k$.

\noindent Definition 2, together with induction on $k$, yields

\textbf{Corollary 2.} $T_{A}(x_{1}^{r_{1}}\cdots x_{k}^{r_{k}})=\underset
{(c_{i,j})\in C(r_{1},\cdots,r_{k})}{\Sigma}\underset{j}{\Pi}(\underset
{i}{\Sigma}c_{i,j})!\underset{i,j}{\Pi}\frac{a_{i,j}^{c_{i,j}}}{c_{i,j}!}$.

\

\textbf{2.3. The formula.} It is well known that simply connected semi-simple
Lie groups are classified by their Cartan matrices [Hu, p.55]. So,
conceivably, any geometric invariant associated to $G/H$ \ can be reduced in
principle to Cartan numbers of $G$ (entries in the Cartan matrix of $G$).
Explicit and direct relationship may become more desirable if one wants to
find an expression of the invariant in its natural generality (i.e. uniformly
for all $G/H$) rather than for special cases. We present both a formula and an
algorithm, which evaluate the number $a_{u,v}^{w}$ in term of Cartan numbers
of $G$.

\bigskip

Assume that $w=\sigma_{\beta_{1}}\circ\cdots\circ\sigma_{\beta_{k}}$,
$\beta_{i}\in\Delta$, is a reduced decomposition of an $w\in\overline{W}$, and
let $A_{w}=(a_{i,j})_{k\times k}$ be the associated Cartan matrix. For a
subset $L=[i_{1},\cdots,i_{r}]\subseteq\lbrack1,\cdots,k]$ we put $\mid
L\mid=r$ and set

$\qquad\sigma_{L}=\sigma_{\beta_{i_{1}}}\circ\cdots\circ\sigma_{\beta_{i_{r}}%
}\in W$;$\quad x_{L}=x_{i_{1}}\cdots x_{i_{r}}\in\mathbb{Z}[x_{1},\cdots
,x_{k}]$.

\noindent Our solution to the problem is

\textbf{Theorem.} If $u,v\in\overline{W}$ with $l(w)=l(u)+l(v)$, then

$\qquad a_{u,v}^{w}=T_{A_{w}}[(\sum\limits_{\substack{\mid L\mid
=l(u)\\\sigma_{L}=u}}x_{L})(\sum\limits_{\substack{\mid K\mid=l(v)\\\sigma
_{K}=v}}x_{K})]$,

\noindent where $L,K\subseteq\lbrack1,\cdots,k]$.

\bigskip

The proof of the Theorem will be developed in such a way as to suggest its
analogue for multiplying generalized Schubert classes in the focal manifolds
of isoparametric submanifolds. In particular, the Theorem is valid for the
focal manifolds of an isoparametric submanifold with equal multiplicities $2$,
in which the classical flag manifolds $G/H$ are special cases (cf.
\textbf{7.2}, \textbf{7.5 }and [HPT]).

In [Bi] S. Billey obtained a recurrence [Bi, (5.5)] that can be used to derive
an expression for $a_{u,v}^{w}$ by using all of the quantities $\xi^{v}%
(t)\mid_{\alpha}$, $\xi^{s}(t)\mid_{\alpha}$ and $\pi_{t}\mid_{\alpha}$ with
$u\leq s<t\leq w$ , where the $\xi^{s}(t)$ and the $\pi_{t}$ for $t,s\in W$
are certain polynomials in the simple roots of $G$ defined respectively by
Kostant and Kumar in [KK] and by Billey in [Bi], and where $\mid_{\alpha}$
means evaluating the polynomials at the regular point $\alpha$. In view of
[Bi, Theorem 4], our theorem expresses $a_{u,v}^{w}$ only in terms of the data
required to describe $\xi^{u}(w)$ and $\xi^{v}(w)$.

\bigskip

\textbf{Example 3.} Let $W$ be the Weyl group of $G_{2}$. Continuing from
Example 1 we express the product $P_{u}P_{v}$ in $H^{\ast}(G_{2}/T)$ in terms
of Schubert classes, where $u$, $v\in W$ are specified by the vectors
$u(\alpha)$, $v(\alpha)\in L(T)$ in the figure. We note that $l(u)=3$,
$l(v)=2$ and $w$, $w^{\prime}\in W$ are the only elements with length
$l(u)+l(v)=5$.

Referring to the reduced decomposition $w=\sigma_{\beta_{2}}\circ\sigma
_{\beta_{1}}\circ\sigma_{\beta_{2}}\circ\sigma_{\beta_{1}}\circ\sigma
_{\beta_{2}}$ obtained in Example 1, the solutions in $L\subset\lbrack
1,2,3,4,5]$ to the equations $\sigma_{L}=u$, $\mid L\mid=3$ are

$\qquad L=(1,2,3),(1,2,5),(1,4,5),(3,4,5)$;

\noindent and the solutions in $K\subset\lbrack1,2,3,4,5]$ to the equations
$\sigma_{K}=v$, $\mid K\mid=2$ are

$\qquad K=(2,3),(2,5),(4,5)$.

\noindent Using the Theorem we compute

$a_{u,v}^{w}=T_{A_{w}}[(x_{1}x_{2}x_{3}+x_{1}x_{2}x_{5}+x_{1}x_{4}x_{5}%
+x_{3}x_{4}x_{5})(x_{2}x_{3}+x_{2}x_{5}+x_{4}x_{5})]$

\qquad$=2+2T_{A_{w}}(x_{1}x_{2}x_{4}x_{5}^{2})+T_{A_{w}}(x_{1}x_{4}^{2}%
x_{5}^{2})$

\qquad$\quad+T_{A_{w}}(x_{2}x_{3}^{2}x_{4}x_{5})+T_{A_{w}}(x_{2}x_{3}%
x_{4}x_{5}^{2})+T_{A_{w}}(x_{3}x_{4}^{2}x_{5}^{2})$,

\noindent where the second equality follows from the additivity of $T_{A_{w}}$
and an application of Corollary 1. With the matrix $A_{w}$ being determined in
Example 1, we find that

\noindent$%
\begin{array}
[c]{ccccc}%
T_{A_{w}}(x_{1}x_{2}x_{4}x_{5}^{2}) & T_{A_{w}}(x_{1}x_{4}^{2}x_{5}^{2}) &
T_{A_{w}}(x_{2}x_{3}^{2}x_{4}x_{5}) & T_{A_{w}}(x_{2}x_{3}x_{4}x_{5}^{2}) &
T_{A_{w}}(x_{3}x_{4}^{2}x_{5}^{2})\\
1 & -2 & 1 & -1 & -1
\end{array}
$.

\noindent Consequently, $a_{u,v}^{w}=1$.

Similarly, from the reduced decomposition $w\prime=\sigma_{\beta_{1}}%
\circ\sigma_{\beta_{2}}\circ\sigma_{\beta_{1}}\circ\sigma_{\beta_{2}}%
\circ\sigma_{\beta_{1}}$ we find

$\sum\limits_{\substack{\mid L\mid=l(u)\\\sigma_{L}(\alpha)=u(\alpha)}%
}x_{L}=x_{2}x_{3}x_{4},\qquad\sum\limits_{\substack{\mid K\mid=l(v)\\\sigma
_{K}(\alpha)=v(\alpha)}}x_{K}=x_{1}x_{2}+x_{1}x_{4}+x_{3}x_{4}$.

\noindent From the Theorem we get

$\qquad a_{u,v}^{w^{\prime}}=T_{A_{w}}[(x_{2}x_{3}x_{4})(x_{1}x_{2}+x_{1}%
x_{4}+x_{3}x_{4})]=0$ (by Corollary 1).

Summarizing, $P_{u}P_{v}=P_{w}$.

\bigskip

\textbf{Remark 1.} It follows from intersection theory that the coefficients
$a_{u,v}^{w}$ are always non-negative. It is therefore attempted to have a
method to compute the number without cancellation involved (i.e.
\textsl{positively multiplying Schubert classes}). From the computation in
Example 3 one finds that our method is not positive. However, the geometric
reason behind this phenomenon can be easily clarified. This will be discussed
in subsection 7.6.

\bigskip

\textbf{2.4. The algorithm.} In concrete situations one prefers to see the
practical value of $a_{u,v}^{w}$ rather than the closed formula, for this
could reveal in a direct way the intersection multiplicities of $X_{u}$ with
$X_{v}$ in the variety $X_{w}$. For this purpose the Theorem does indicate an
effective algorithm to evaluate $a_{u,v}^{w}$, as the following recipe shows
(see also Examples 1-3).

(1) starting from the Cartan matrix of $G$, \ a program to enumerate all
elements in a coset $\overline{W}$ of the Weyl group $W$ by their
\textsl{minimal reduced decompositions} is available in [DZZ];

(2) for an $w\in\overline{W}$ with a reduced decomposition, the corresponding
Cartan matrix $A_{w}$ can be read directly from Cartan matrix of $G$ (cf.
Example 1);

(3) for an $w\in\overline{W}$ with a reduced decomposition $w=\sigma
_{\beta_{1}}\circ\cdots\circ\sigma_{\beta_{k}}$ and an $u\in\overline{W}$ with
$l(u)=r<k$, the solutions in the subsequence $[i_{1},\cdots,i_{r}%
]\subseteq\lbrack1,\cdots,k]$ to the equation $\sigma_{\beta_{i_{1}}}%
\circ\cdots\circ\sigma_{\beta_{i_{r}}}=u$ in $W$ agree with the solutions to
the vector equation $\sigma_{\beta_{i_{1}}}\circ\cdots\circ\sigma
_{\beta_{i_{r}}}(\alpha)=u(\alpha)$ in the linear space $L(T)$, where
$\alpha\in L(T)$ is a fixed regular point;

(4) the evaluation the operator $T_{A_{w}}$ on a polynomial can be easily
programmed (cf. Example 2).

Combining the ideas explained above, a program to compute the numbers
$a_{u,v}^{w}$ has been compiled [DZ$_{2}$]. It uses Cartan matrix as the only
input and computations in various flag manifolds $G/H$ can be performed by the
single program.

As for the efficiency of the program, we refer the reader to the computational
results tabulated in [DZ$_{1}$]. They were produced by a similar program that
implements Steenrod operations on Schubert classes.

\bigskip

\textbf{2.5. Arrangement of the paper.} The rest sections of the paper are so
arranged. Section 3 develops preliminary results from algebraic topology. We
recall from [Du$_{2}$] the cohomology of twisted product of $2$-spheres (Lemma
3.3), and the enumerative formula on these manifolds (Lemma 3.4). In
particular, we introduce \textsl{divided differences} for spherical
represented involutions, and their basic properties are established in Lemma 3.2.

We shall see in Section 4 that, by resorting to the geometry of the adjoint
representation, Bott-Samelson cycles in the space $G/T$ appears as certain
twisted product of $2$-spheres that are parameterized by ordered sequences of
roots, and the divided differences on the integral cohomology of $G/T$ arise
naturally from the geometric fact that the involution on $G/T$ corresponding
to a root is spherical representable. After determining the induced action of
Bott-Samelson cycles (corresponding to a sequence of simply roots) on Schubert
classes in Lemma 5.1 (Section 5), the Theorem is established in Section 6.

\bigskip

In many literatures ranging from topology, algebraic and differential geometry
to representation theory, one finds analogues of the terminologies that we
work with, such as Bott-Samelson cycles (or schemes), divided differences and
Schubert varieties, but with seemingly different appearances. In particular,
Schubert varieties were originally introduced and extensively studied in the
context of algebraic geometry, but we will work with Lie groups in the real
compact form so that our method are ready to apply to general situations (cf.
\textbf{7.5}). In order to merge our presentations into the existing
literatures and, at the same time, not to interrupt our exposition, Section 7,
entitled \textsl{Historical remarks}, is devoted to recall those historical
events illustrating the readiness and necessity of the conceptual development
in our paper, and will be referred to from time to time.

\bigskip

Finally, a brief account for the method of the proof. In 1973 Hansen
discovered that the celebrated K-cycles on the flag manifold $G/T$ constructed
by Bott and Samelson in 1958 (cf. [BS$_{2}$] or \textbf{7.1}) provided a
degree $1$ map $g_{w}$ from a twisted products $\Gamma_{w}$ of $2$-spheres
onto the Schubert variety $X_{w}$ (cf. [Han] or \textbf{7.4}). This suggests
that the intersection product in $X_{w}$ can be translated as part of the
intersection product in $\Gamma_{w}$ via the homomorphism induced by $g_{w}$.
However, the latter is much easier to work with for the following reasons (cf.
Lemma 3.3 and Lemma 3.4 in Section 3).

a) the space $\Gamma_{w}$ admits a natural cell decomposition\textsl{\ }with
each cell, again,\textsl{\ a twisted product of $2$-spheres};

b) the cohomology of $\Gamma_{w}$ is a polynomial ring $\mathbb{Z}%
[x_{1},\cdots,x_{k}]$ generated by $x_{i}$'s in dimension $2$, subject to
relations\textsl{\ }occurring\textsl{ only in dimension }$4$;

c) the intersection product in $\Gamma_{w}$ is handled by a triangular
operator $T_{A}$.

\noindent Therefore, the intersection multiplicity $a_{u,v}^{w}$ in question
can be calculated by computations in a space like $\Gamma_{w}$, rather than in
the Schubert variety $X_{w}$ itself.

\section{Preliminaries in topology}

In this paper all homologies (resp. cohomologies) will have integer
coefficients unless otherwise stated. If $f:X\rightarrow Y$ is a continuous
map between two topological spaces, $f_{\ast}$ (resp. $f^{\ast}$) is the
homology (resp. cohomology) map induced by $f$. Write $S^{r}$ for the
$r$-dimensional sphere. If $M$ is an oriented closed manifold (resp. a
connected projective variety) $[M]\in H_{\dim M}(M)$ stands for the
orientation class. The Kronecker pairing, between cohomology and homology of a
space $X$, will be denoted by $<,>:H^{\ast}(X)\times H_{\ast}(X)\rightarrow
\mathbb{Z}$.

\bigskip

\textbf{3.1. Sphere bundle with a cross section. }Let $p:E\rightarrow M$ be a
smooth, oriented $r$-sphere bundle over an oriented manifold $M$ which has a
cross section $s:M\rightarrow E$. Let the normal bundle $\xi$ of the embedding
$s$ be oriented by $p$, and let $e\in H^{r}(M)$ be the Euler class of $\xi$
with respect to this orientation.

The integral cohomology $H^{\ast}(E)$ can be described as follows. Denote by
$i:S^{r}\rightarrow E$ for the fiber inclusion of $p$ over a point $z\in M$,
and write by $J:$ $E\rightarrow E$ for the involution given by the antipodal
map in each fiber sphere. We have the following result from [Du$_{2}$, Lemma 4].

\textbf{Lemma 3.1. }\textsl{There exists a unique class }$x\in H^{r}%
(E)$\textsl{ such that}

\textsl{\qquad}$s^{\ast}(x)=0\in H^{\ast}(M)$\textsl{ and }$<i^{\ast
}(x),[S^{r}]>=1$\textsl{. }

\noindent\textsl{Furthermore }

\textsl{(1) }$H^{\ast}(E)$\textsl{, as a module over }$H^{\ast}(M)$\textsl{,
has the basis }$\{1,x\}$\textsl{ subject to the relation }$x^{2}+p^{\ast
}(e)x=0$\textsl{;}

\textsl{(2) the induced cohomology map }$J^{\ast}$\textsl{ acts}
\textsl{identically on the subset} $\operatorname{Im}p^{\ast}\subset H^{\ast
}(E)$ \textsl{and}

\textsl{\qquad\qquad}$J^{\ast}(x)=(-1)^{r-1}x-p^{\ast}(e)$\textsl{.}

\textbf{Remark 2.} If $r$ is odd, $2e=0$.

\textsl{\bigskip}

\textbf{3.2. Divided difference of a spherical represented involution. }A
self-map $\sigma$ of a manifold $M$ is called \textsl{an involution} if
$\sigma^{2}=id:M\rightarrow M$. An $r$-\textsl{spherical representation} of
the involution $(M;\sigma)$ is a system $f:(E;J)\rightarrow(M;\sigma)$ in which

(1) $E$ is the total space of an oriented $r$-sphere bundle $p:E\rightarrow M$
with a cross section $s$;

(2) $f$ is a continuous map $E\rightarrow M$ that satisfies the following two constrains

\noindent(3.1)$\qquad f\circ s=id:M\rightarrow M$; and

\noindent(3.2)$\qquad f\circ J=\sigma\circ f:$ $E\rightarrow M$,

\noindent where $J:$ $E\rightarrow E$ is the involution on $E$ given by the
antipodal map in the fibers.

\bigskip

In view of the $H^{\ast}(M)$-module structure of $H^{\ast}(E)$ specified by
Lemma 3.1, an $r$-spherical representation $f$ of the involution $(M,\sigma)$
gives rise to an additive operator $\theta_{f}:H^{m}(M)\rightarrow H^{m-r}(M)$
of degree $-r$ that is characterized uniquely as follows.

\textsl{The induced homomorphism }$f^{\ast}:H^{\ast}(M)\rightarrow H^{\ast
}(E)$\textsl{ satisfies}

\noindent(3.3)$\qquad\qquad\qquad f^{\ast}(z)=p^{\ast}(z)+p^{\ast}(\theta_{f}(z))x$

\noindent\textsl{for all }$z\in H^{\ast}(M)$\textsl{.}

\bigskip

Useful properties of $\theta_{f}$ can be derived directly from its definition (3.3).

\textbf{Lemma 3.2.} \textsl{Let }$e\in H^{r}(M)$\textsl{ be as that in Lemma
3.1. We have}

\textsl{(1) }$\sigma^{\ast}=Id-e\theta_{f}$\textsl{ }$:$\textsl{ }$H^{\ast
}(M)\rightarrow H^{\ast}(M)$\textsl{;}

\textsl{(2) }$\theta_{f}(z_{1}z_{2})=\theta_{f}(z_{1})z_{2}+\sigma^{\ast
}(z_{1})\theta_{f}(z_{2})$\textsl{, }$z_{1},z_{2}\in H^{\ast}(M)$\textsl{.}

\textsl{(3) if }$r$\textsl{ is even, then }$\theta_{f}(e)=2$\textsl{(}$\in
H^{0}(M)=\mathbb{Z}$\textsl{); if }$r$\textsl{ is odd, then }$\theta_{f}(e)=0$\textsl{;}

\textsl{(4) if }$r$\textsl{ is even, }$2\theta_{f}\circ\theta_{f}=0:$\textsl{
}$H^{\ast}(M)\rightarrow H^{\ast}(M)$\textsl{.}

\textbf{Proof.} Since $f^{\ast}$ is a ring map, we have

$f^{\ast}(z_{1}z_{2})=[p^{\ast}(z_{1})+p^{\ast}(\theta_{f}(z_{1}))x][p^{\ast
}(z_{2})+p^{\ast}(\theta_{f}(z_{2}))x]$

\qquad$=p^{\ast}(z_{1}z_{2})+p^{\ast}[\theta_{f}(z_{1})z_{2}+z_{1}\theta
_{f}(z_{2})]x+p^{\ast}[\theta_{f}(z_{1})\theta_{f}(z_{2})]x^{2}$

\qquad$=p^{\ast}(z_{1}z_{2})+p^{\ast}[\theta_{f}(z_{1})z_{2}+z_{1}\theta
_{f}(z_{2})-e\theta_{f}(z_{1})\theta_{f}(z_{2})]x$,

\noindent where the last equality is obtained from $x^{2}=-p^{\ast}(e)x$
(because of the relation $x^{2}+p^{\ast}(e)x=0$). Comparing this with (3.3) yields

\noindent(3.4)$\qquad\qquad\theta_{f}(z_{1}z_{2})=\theta_{f}(z_{1})z_{2}%
+z_{1}\theta_{f}(z_{2})-e\theta_{f}(z_{1})\theta_{f}(z_{2})$

\qquad\qquad\qquad\qquad\ \noindent$=\theta_{f}(z_{1})z_{2}+(z_{1}-e\theta
_{f}(z_{1}))\theta_{f}(z_{2})$.

Applying $J^{\ast}$ to (3.3) gives

$\qquad J^{\ast}f^{\ast}(z)=p^{\ast}(z)+p^{\ast}(\theta_{f}(z))J^{\ast}(x)$

\noindent by (2) of Lemma 3.1. From $J^{\ast}f^{\ast}=f^{\ast}\sigma^{\ast}$
(by (3.2)) and $J^{\ast}(x)=(-1)^{r-1}x-p^{\ast}(e)$ (by (2) of Lemma 3.1) we get

$\qquad f^{\ast}(\sigma^{\ast}(z))=p^{\ast}(z-e\theta_{f}(z))+p^{\ast
}((-1)^{r-1}\theta_{f}(z))x$.

\noindent Comparing this with (3.3) gives rise to

\noindent(3.5)$\qquad\sigma^{\ast}(z)=z-e\theta_{f}(z)$; and

\noindent(3.6)$\qquad\theta_{f}(\sigma^{\ast}(z))=(-1)^{r-1}\theta_{f}(z)$.

\noindent(1) is verified by (3.5). Combining (3.4) with (3.5) shows (2).

Substituting (3.5) in the left hand side of (3.6) and rewriting the resulting
equation by using (3.4) gives

\noindent(3.7)$\qquad\qquad(1+(-1)^{r})\theta_{f}(z)=\theta_{f}(e)\theta
_{f}(z)+(e-\theta_{f}(e)e)\theta_{f}(\theta_{f}(z))$.

\noindent Taking $z=e$ in (3.7) and noting that $\theta_{f}(\theta_{f}(e))\in
H^{-r}(M)=0$ we get the equation $(1+(-1)^{r})\theta_{f}(e)=\theta_{f}(e)^{2}$
in $H^{0}(M)=\mathbb{Z}$. This proves (3).

Now (3.7) becomes $e\theta_{f}(\theta_{f}(z))=0$ for all $z\in H^{\ast}(M)$ by
(3). Consequently, taking $z=\theta_{f}(u)$ in (3.5) yields $\sigma^{\ast
}(\theta_{f}(u))=\theta_{f}(u)$. Now (3.6) implies that if $r$ is even,
$2\theta_{f}(\theta_{f}(u))=0$ for all $u\in H^{\ast}(M)$. This verifies (4),
hence completes the proof of Lemma 3.2.

\bigskip

We observe from (1) of Lemma 3.2 that, for any $z\in H^{\ast}(M)$, the
difference $z-\sigma^{\ast}(z)\in H^{\ast}(M)$ is always divisible by $e$ with
quotient $\theta_{f}(z)$.

\textbf{Definition 3}. The operator $\theta_{f}$ is called the \textsl{divided
difference } of the spherical representation $f$ of the involution
$(M,\sigma)$ (Compare with the discussion in\textbf{ 7.5}).

\bigskip

\textbf{3.3. Twisted product of 2-spheres. }The following definition singles
out a class spaces in which we will be particularly interested.\textsl{\ }

\textbf{Definition 4.} A smooth manifold $M$ is called an \textsl{oriented
twisted product of 2-spheres of rank k}, denoted by $M=\underset{1\leq i\leq
k}{\propto}S^{2}$, if there is a tower of smooth maps

$\qquad\qquad M=M_{k}\overset{p_{k-1}}{\rightarrow}M_{k-1}\overset{p_{k-2}%
}{\rightarrow}\cdots\overset{p_{2}}{\rightarrow}M_{2}\overset{p_{1}%
}{\rightarrow}M_{1}$

\noindent in which

1) $M_{1}$ is diffeomorphic to $S^{2}$ with an orientation;

2) $p_{i}$ is the projection of an oriented smooth $S^{2}$ bundle over $M_{i}
$;

3) $p_{i}$ has a fixed cross section $s_{i}$, $1\leq i\leq k-1$.

\smallskip\

Let $M=\underset{1\leq i\leq k}{\propto}S^{2}$ be a twisted product of
$2$-spheres of rank $k$. The cross sections $s_{i}$ yield a sequence of embeddings

$\qquad\qquad M_{1}\overset{s_{1}}{\rightarrow}M_{2}\overset{s_{2}%
}{\rightarrow}\cdots\overset{s_{k-1}}{\rightarrow}M_{k}=M$.

\noindent In view of this we shall make no notational distinction between a
subspace $N\subset M_{i}$ and its image under $s_{i}$ in $M_{j}$, $j\geq i$.
Take a base point $x_{0}\in$ $M_{1}$ and let each $M_{i}$ has $x_{0}$ as its
base point.

For a subset $L=$ $[i_{1},\cdots,i_{r}]\subseteq\lbrack1,\cdots,k]$ a smooth submanifold

$\qquad\quad\quad S(L)\subset M_{i_{k}}\subset M_{i_{k}+1}\subset\cdots\subset M$

\noindent can be introduced inductively as follows.

1) $S(1)=M_{1}$;

2) if $i>1$, $S(i)\subset M_{i}$ is the fiber sphere of $p_{i-1}$ over the
base point;

3) assume that $S(L^{^{\prime}})\subset M_{i_{r-1}}$, where $L^{\prime}%
=[i_{1},\cdots,i_{r-1}]$, has been defined and consider the case $L=$
$[i_{1},\cdots,i_{r-1},j]$. Then $S(L)\subset M_{j}$ is the total space of the
restricted bundle of $p_{j-1}:M_{j}\rightarrow M_{j-1}$ to the subspace
$S(L^{^{\prime}})\subset M_{j-1}$. The natural bundle map over the inclusion
$S(L^{^{\prime}})\subset M_{i_{r-1}}\rightarrow M_{j-1}$ gives rise to the
desired embedding $S(L)\subseteq M_{j}\subseteq M$.

\bigskip

The integral homology (resp. cohomology ) of an $M=\underset{1\leq i\leq
k}{\propto}S^{2}$ can be described as follows. Consider the normal bundle
$\pi_{i}:E_{i}\rightarrow M_{i-1}$ of the embedding $s_{i-1}:M_{i-1}%
\rightarrow M_{i}$, $2\leq i\leq k$. It is a $2$-plane bundle with a natural
orientation inherited from that on $p_{i-1}$. Let $e_{i}^{\prime}\in
H^{2}(M_{i-1})$ be the Euler class of $\pi_{i}$. We set $e_{1}=0$ and for
$i\geq2$

$\qquad\qquad e_{i}=(p_{i-1}\circ p_{i}\circ\cdots\circ p_{k-1})^{\ast}%
e_{i}^{\prime}\in H^{2}(M)$.

\bigskip

\textbf{Lemma 3.3}. \textsl{For an }$M=\underset{1\leq i\leq k}{\propto}S^{2}%
$\textsl{ let }$[S(L)]\in H_{2\mid L\mid}(M)$\textsl{ be the fundamental class
of the cycle }$S(L)\subset M$\textsl{. Then}

\textsl{(1) the set }$\{[S(L)]\mid L\subseteq\lbrack1,\cdots,k]\}$\textsl{ of
homology classes is an additive basis for the graded }$\mathbb{Z}%
$\textsl{-module }$H_{\ast}(M)$\textsl{.}

\textsl{Further,} \textsl{let} $\{x_{L}\in H^{2\mid L\mid}(M)\mid
L\subseteq\lbrack1,\cdots,k]\}$\textsl{ be the basis of }$H^{\ast}(M)$\textsl{
Kronecker dual to the basis }$\{[S(L)]\mid L\subseteq\lbrack1,\cdots
,k]\}$\textsl{ in homology. Then}

\textsl{(2) }$x_{L}=x_{i_{1}}\cdots x_{i_{r}}$\textsl{ if }$L=$\textsl{
}$[i_{1},\cdots,i_{r}]$\textsl{;}

\textsl{(3) }$H^{\ast}(M)=\mathbb{Z}[x_{1},\cdot\cdot\cdot,x_{k}]/<x_{i}%
^{2}+e_{i}x_{i};$\textsl{ }$1\leq i\leq k>$\textsl{, where }$e_{i}$\textsl{ is
a polynomial in }$x_{1},\cdot\cdot\cdot,x_{i-1}$\textsl{. }

\textbf{Proof.} By setting $T_{r}(M)=\underset{\mid L\mid=r}{\cup
}S(L)\subseteq M$, $1\leq r\leq k$, we get a filtration of subspaces

\noindent(3.8)$\qquad\qquad T_{1}(M)\subset T_{2}(M)\subset\cdots\subset
T_{k}(M)=M$.

\noindent Moreover, the subspace $T_{r}(M)\backslash T_{r-1}(M)$ consisting of
$2r$-dimensional open cells in one-to-one correspondence with the subset
$L\subseteq\lbrack1,\cdots,k]$ with $\mid L\mid=r$. This implies that (3.8)
dominates $M$ by a cell complex with only even dimensional cells. This proves (1).

Assertions (2) and (3) follow easily from (1) of Lemma 3.1 together with
induction on $k$.

\bigskip

For degree reasons the polynomials $e_{i}$ are all homogeneous of degree $1$
in $x_{1},\cdot\cdot\cdot,x_{i-1}$, i.e.

\qquad\qquad$e_{i}=a_{1,i}x_{1}+\cdots+a_{i-1,i}x_{i-1}$, $a_{i,j}%
\in\mathbb{Z},$ $1\leq i\leq k$.

\textbf{Definition 5.} With $a_{i,j}=0$ for $i\geq j$ being understood, the
strictly upper triangular matrix $A=(-a_{i,j})_{k\times k}$ is called
\textsl{the structure matrix} of $M=\underset{1\leq i\leq k}{\propto}S^{2}$.

\textbf{Remark 3.} It was shown in [Du$_{2}$, Proposition 1] that any strictly
upper triangular matrix $A$ of rank $k$ can be realized as the structure
matrix of an $M=\underset{1\leq i\leq k}{\propto}S^{2}$.

\bigskip

\textbf{3.4. Integration along a twisted product of 2-spheres. }Recall from
Lemma 3.3 that the integral cohomology of an $M=\underset{1\leq i\leq
k}{\propto}S^{2}$ can be specified as a quotient of a free polynomial ring

$\qquad H^{\ast}(M)=\mathbb{Z}[x_{1},\cdot\cdot\cdot,x_{k}]/<x_{i}^{2}%
+e_{i}x_{i};$\textsl{ }$1\leq i\leq k>$.

\noindent Write $\mathbb{Z}[x_{1},\cdot\cdot\cdot,x_{k}]^{(r)}$ for the subset
of all homogeneous polynomials of degree $r$ in $\mathbb{Z}[x_{1},\cdot
\cdot\cdot,x_{k}]$ and let $p_{M}:\mathbb{Z}[x_{1},\cdot\cdot\cdot
,x_{k}]^{(r)}\rightarrow H^{2r}(M)$ be the obvious quotient map. Consider the
additive correspondence $\int_{M}:\mathbb{Z}[x_{1},\cdot\cdot\cdot
,x_{k}]^{(k)}\rightarrow\mathbb{Z}$ defined by

$\qquad\qquad\int_{M}h=<p_{M}(h),[M]>$,

\noindent where $[M]\in H_{2k}(M)=\mathbb{Z}$ is the orientation\textsl{\ }class.

As being indicated by the notation, the operator $\int_{M}$ can be interpreted
as ``\textsl{integration along }$M$'' in De Rham theory. In view of
discussions in [Du$_{1}$], the problem of effective computation of the
intersection product in $M$ asks an effective algorithm to evaluate $\int_{M}%
$. The idea of structure matrix of $M$, together with the operator $T_{A}$
introduced in \textbf{2.2}, is useful in presenting such an algorithm. The
following result was shown in [Du$_{2}$, Proposition 2].

\textbf{Lemma 3.4}. If $M$ has the structure matrix $A=(a_{ij})_{k\times k}$, then

\qquad\qquad\ \ $\int_{M}=T_{A}:\mathbb{Z}[x_{1},\cdot\cdot\cdot,x_{k}%
]^{(k)}\rightarrow\mathbb{Z}$.\

\section{Geometry from adjoint representation}

\smallskip\ Since the point $\alpha\in L(T)$ is regular, the adjoint
representation $Ad:G\rightarrow L(G)$ yields a smooth embedding

$\qquad\varphi:G/T\rightarrow L(G)\quad\quad$by $\varphi(gT)=Ad_{g}(\alpha)$.

\noindent In this way $G/T$ becomes a submanifold of the Euclidean space
$L(G)$.

\bigskip

\ \textbf{4.1.} \textbf{Preliminaries.} Let $\Phi^{+}$ be the set of positive
roots relative to $\alpha$ (cf. \textbf{2.1}). Assume that the Cartan
decomposition of the Lie algebra $L(G)$ relative to $T\subset G$ is

$\qquad\ \qquad L(G)=L(T)\oplus_{\beta\in\Phi^{+}}F_{\beta}$,

\noindent where $F_{\beta}$ is the root space, viewed as a real $2$-plane,
belonging to the root $\beta\in\Phi^{+}$ ([Hu, p.35]). Let $[,]$ be the Lie
bracket on $L(G)$. We make a simultaneous choice of an orientation on each
$F_{\beta}$ by the following rule: take a $v\in F_{\beta}\backslash\{0\}$ and
let $v\prime\in F_{\beta}$ be such that $[v,v^{\prime}]=\beta$. Then
$F_{\beta}$ is oriented by the ordered basis $\{v,v^{\prime}\}$. Clearly this
orientation on $F_{\beta}$ is independent of a specific choice of $v$.

\smallskip\

We quote from [HPT, p.426-427] for relevant information concerning the
geometries of the submanifold $G/T\subset L(G)$.

(4.1) The subspaces $\oplus_{\beta\in\Phi^{+}}F_{\beta}$ and $L(T)$ of $L(G)$
are tangent and normal to $G/T$ at $\alpha$ respectively;

(4.2) The tangent bundle to $G/T$ has a canonical orthogonal decomposition
into the sum of $m$ integrable $2$-plane bundles $\oplus_{\beta\in\Phi^{+}%
}E_{\beta}$ with $E_{\beta}(\alpha)=F_{\beta}$.

In general, $E_{\beta}(y)=Ad_{g}(F_{\beta})$ if $y=Ad_{g}(\alpha)\in G/T$.

(4.3) The leaf of the integrable subbundle $E_{\beta}$ through a point $y\in
G/T$, denoted by $S(y;\beta)$, is a $2$-sphere.

\noindent In view of (4.2) and (4.3) we adopt the following convention.

\textbf{Convention.} Each bundle $E_{\beta}$, $\beta\in\Phi^{+}$, is oriented
such that the identification $E_{\beta}(\alpha)=F_{\beta}$ is orientation
preserving. Let the $2$-sphere $S(y;\beta)$ have the orientation inherited
from that on its tangent plane $E_{\beta}(y)$ at $y\in S(y;\beta)$.

\bigskip

By considering $G/T$ as a submanifold in the Euclidean space $L(G)$ via the
embedding $\varphi$, the geometric features singled out in (4.2) and (4.3)
will dominate all of our geometric constructions and computations in this section.

\bigskip

\textbf{4.2.} \textbf{Realization of roots as homology and cohomology classes.
}We start by relating the set $\Phi$ of roots with certain $2$-dimensional
homology (resp. cohomology) classes of $G/T$.

For a $\beta\in\Phi^{+}$ write by $e_{\beta}\in H^{2}(G/T)$ for the Euler
class of the oriented bundle $E_{\beta}$, and let $\varphi_{\beta}%
:S(y;\beta)\rightarrow G/T$ be inclusion of the leaf sphere (cf. (4.3)).

\textbf{Lemma 4.1.} \textsl{The Kronecker paring }$H^{2}(G/T)\times
H_{2}(G/T)\rightarrow\mathbb{Z}$\textsl{ can be expressed in term of Cartan
numbers as}

$\qquad\qquad<e_{\gamma},\varphi_{\beta\ast}[S(y;\beta)]>=\beta\circ\gamma
$\textsl{, }$\beta,\gamma\in\Phi$\textsl{,}

\noindent\textsl{where }$[S(y;\beta)]\in H^{2}(S(y;\beta))$\textsl{(}%
$=\mathbb{Z}$\textsl{) is the orientation class.}

\textbf{Remark 4.} The class $\varphi_{\beta\ast}[S(y;\beta)]\in H_{2}(G/T)$
is independent of the choice of the point $y\in G/T$. Let $\varepsilon
:[0,1]\rightarrow G/T$ be a path that joins $y$ to some $y^{\prime}$. Then the
continuous one-parameter family of $2$-spheres $S(\varepsilon(t);\beta)\subset
G/T$, $t\in\lbrack0,1]$ is an isotopy from $S(y;\beta)$ to $S(y^{\prime}%
;\beta)$.

\bigskip

The Weyl group $W$ of $G$, acting as isometries of the Cartan subalgebra
$L(T)$, has the effect to permute the set $\Phi$ of roots [Hu]. On the other
hand, via the embedding $\varphi$, the canonical action of $W$ on $G/T$
[BS$_{2}$] is given by

\noindent(4.4)$\qquad\qquad w(x)=Ad_{g}(w(\alpha))$, $x=Ad_{g}(\alpha)\in
G/T$, $w\in W$.

\noindent Let $w^{\ast}:H^{\ast}(G/T)\rightarrow H^{\ast}(G/T)$ be the induced
action on cohomology.

\textbf{Lemma 4.2.} \textsl{For an }$w\in W$\textsl{ one has}

\textsl{(1)} $w(S(y;\beta))=S(w(y);w(\beta))$\textsl{;}

\textsl{(2) }$w^{\ast}(E_{\beta})=E_{w^{-1}(\beta)}$\textsl{; and}

\textsl{(3)} $w^{\ast}(e_{\beta})=e_{w^{-1}(\beta)}$\textsl{.}

\textbf{Proof.} In term of the $W$-action on the set $\Phi$ of roots, the
induced bundle $w^{\ast}(E_{\beta})$ is $E_{w^{-1}(\beta)}$ (cf. [HTP,
\textbf{1.6}]). This shows (2). Assertion (3) comes now from the naturality of
Euler classes.

Item (1), which follows also from (2) (and (4.3)), indicates that the action
of $w$ on $G/T$ has the effect to carry the leaf sphere through $y$
corresponding to a root $\beta$ diffeomorphically onto the leaf sphere through
$w(y)$ corresponding to the root $w(\beta)$.

\bigskip

\textbf{4.3. }$2$-\textbf{spherical represented involutions. }Each
root\textbf{ }$\beta\in\Phi^{+}$ gives rise to an involution $\sigma_{\beta
}:G/T\rightarrow G/T$ in the fashion of (4.4), and defines also the subspace

$\qquad\qquad S(\beta)=\{(y,y_{1})\in G/T\times G/T\mid y_{1}\in S(y;\beta)\}$.

\noindent Projection $p_{\beta}:S(\beta)\rightarrow G/T$ onto the first factor
is easily seen to be a $2$-sphere bundle projection (with the leaf sphere
$S(y;\beta)$ as the fiber over $y\in G/T$). The map $s_{\beta}:G/T\rightarrow
S(\beta)$ by $s_{\beta}(y)=(y,y)$ furnishes $p_{\beta}$ with a ready-made
cross section.

Let $J_{\beta}$ be the involution on $S(\beta)$ given by the antipodal map on
each fiber sphere and let $f_{\beta}:S(\beta)\rightarrow G/T$ be the
projection onto the second factor. Then, as is clear,

$\qquad\qquad f_{\beta}\circ s_{\beta}=id:G/T\rightarrow G/T$;

$\qquad\qquad f_{\beta}\circ J_{\beta}=\sigma_{\beta}\circ f_{\beta}%
:S(\beta)\rightarrow G/T.$

\noindent That is, \textsl{the map }$f_{\beta}:(S(\beta),J_{\beta}%
)\rightarrow(G/T,\sigma_{\beta})$\textsl{ is a }$2$\textsl{-spherical
representation of the involution }$(G/T,\sigma_{\beta})$\textsl{ }(cf.
\textbf{3.2}).

\bigskip

\textbf{Lemma 4.3.} \textsl{Assume the same setting as the above.}

\textsl{1) As a }$H^{\ast}(G/T)$\textsl{-module, the cohomology of }$S(\beta
)$\textsl{ is given by}

$\qquad\qquad H^{\ast}(S(\beta))=H^{\ast}(G/T)[1,x]/<x^{2}+e_{\beta}x>$\textsl{,}

\noindent\textsl{in which }$x\in H^{2}(S(\beta))$\textsl{ is uniquely
characterized by }$s_{\beta}^{\ast}(x)=0$ \textsl{and}

$\qquad<i^{\ast}(x),[S(y;\beta)]>=1$\textsl{ (for all }$y\in G/T$\textsl{), }

\noindent\textsl{where }$i:S(y;\beta)\rightarrow S(\beta)$\textsl{ is the
inclusion of the fiber over} $y\in G/T$.

\textsl{2) Write }$\theta_{\beta}:H^{\ast}(G/T)\rightarrow H^{\ast}%
(G/T)$\textsl{ to denote the divided difference associated to the 2-spherical
representation }$f_{\beta}$ \textsl{of the involution }$(G/T,\sigma_{\beta}%
)$\textsl{. Then}

\textsl{(1) for any }$\gamma\in\Phi$, $\theta_{\beta}(e_{\gamma})=\beta
\circ\gamma$\textsl{(}$\in H^{0}(G/T)=\mathbb{Z}$\textsl{);}

\textsl{(2) }$\theta_{\beta}\circ\theta_{\beta}=0$\textsl{;}

\textsl{(3) }$\sigma_{\beta}^{\ast}=Id-e_{\beta}\theta_{\beta}$\textsl{ }%
$:$\textsl{ }$H^{\ast}(G/T)\rightarrow H^{\ast}(G/T)$\textsl{;}

\textsl{(4) }$\theta_{\beta}(z_{1}z_{2})=\theta_{\beta}(z_{1})z_{2}%
+\sigma_{\beta}^{\ast}(z_{1})\theta_{\beta}(z_{2})$\textsl{, }$z_{1},z_{2}\in
H^{\ast}(G/T)$\textsl{, and}

\textsl{(5) for any }$w\in W$, $\theta_{w(\beta)}=(w^{-1})^{\ast}\theta
_{\beta}w^{\ast}$.

\textbf{Proof.} Since the normal bundle of the embedding $s_{\beta}$ is
$E_{\beta}$, one has the relation $x^{2}+e_{\beta}x=0$ by (1) of Lemma 3.1.
This verifies 1).

In 2) properties (3) and (4) corresponds to the items (3) and (4) in Lemma
3.2. (2) follows from (2) of Lemma 3.2 since $H^{\ast}(G/T)$ is torsion free
[BS$_{1}$]. It remains to show (1) and (5).

Property (1) is verified by

$\qquad\beta\circ\gamma=<e_{\gamma},\varphi_{\beta\ast}[S(y;\beta)]>$(by Lemma 4.1)

$=<\varphi_{\beta}^{\ast}(e_{\gamma}),[S(y;\beta)]>$ (by the naturality of
$<,>$)

$=<i^{\ast}\circ f_{\beta}^{\ast}(e_{\gamma}),[S(y;\beta)]>$ (since
$\varphi_{\beta}=f_{\beta}\circ i$)

$=<i^{\ast}(p_{\beta}^{\ast}(e_{\gamma})+p_{\beta}^{\ast}(\theta_{\beta
}(e_{\gamma}))x),[S(y;\beta)]>$ (the definition of $\theta_{\beta}$ in (3.3))

$=\theta_{\beta}(e_{\gamma})$ (by 1) of this Lemma),

\noindent where in the last equality we have applied the standard fact that

\noindent(4.5) \textsl{the composition }$i^{\ast}\circ p_{\beta}^{\ast}%
:H^{r}(G/T)\rightarrow H^{r}(S(y;\beta))$\textsl{ is zero in degree }%
$r>0$\textsl{, and is an isomorphism }$\mathbb{Z}\rightarrow\mathbb{Z}%
$\textsl{ if }$r=0$\textsl{.}

Finally we show (5). For an $w\in W$ the self diffeomorphism $w\times
w:G/T\times G/T\rightarrow G/T\times G/T$ restricts to a diffeomorphism
$\widetilde{w}:S(\beta)\rightarrow S(w(\beta))$ (by (1) of Lemma 4.2) that
fits into the following two commutative diagrams

$\qquad%
\begin{array}
[c]{ccc}%
S(\beta) & \overset{\widetilde{w}}{\rightarrow} & S(w(\beta))\\
p_{\beta}\downarrow &  & \downarrow p_{w(\beta)}\\
G/T & \overset{w}{\rightarrow} & G/T
\end{array}
$ and $%
\begin{array}
[c]{ccc}%
S(\beta) & \overset{\widetilde{w}}{\rightarrow} & S(w(\beta))\\
f_{\beta}\downarrow &  & \downarrow f_{w(\beta)}\\
G/T & \overset{w}{\rightarrow} & G/T\text{.}%
\end{array}
$

\noindent From the definition of $\theta_{\beta}$ in (3.3) we have

$\qquad f_{\beta}^{\ast}(w^{\ast}(z))=p_{\beta}^{\ast}(w^{\ast}(z))+p_{\beta
}^{\ast}(\theta_{\beta}(w^{\ast}(z)))x$, $z\in H^{\ast}(G/T)$.

\noindent On the other hand

$f_{\beta}^{\ast}(w^{\ast}(z))=\widetilde{w}^{\ast}f_{w(\beta)}^{\ast}(z)$ (by
the commutivity of the second diagram)

$=\widetilde{w}^{\ast}(p_{w(\beta)}^{\ast}(z)+p_{\beta}^{\ast}[\theta
_{w(\beta)}(z)]x^{\prime})$ (by the definition of $\theta_{w(\beta)}$)

$=p_{\beta}^{\ast}(w^{\ast}(z))+p_{\beta}^{\ast}(w^{\ast}[\theta_{w(\beta
)}(z)])x$ (by the commutivity of the first diagram).

\noindent Comparing the coefficients of $x\in H^{\ast}(S(\beta))$ in the above
two expressions of $f_{\beta}^{\ast}(w^{\ast}(z))$ yields $w^{\ast}%
\theta_{w(\beta)}=\theta_{\beta}w^{\ast}$. This shows (5).

\bigskip

\textbf{Lemma 4.4.} \textsl{Given an ordered sequence }$\beta_{1},\cdots
,\beta_{k}\in\Delta$\textsl{ of simple roots we put }$w=\sigma_{\beta_{1}%
}\circ\cdots\circ\sigma_{\beta_{k}}$\textsl{, }$\theta_{(\beta_{1}%
,\cdots,\beta_{k})}=\theta_{\beta_{1}}\circ\cdots\circ$\textsl{ }%
$\theta_{\beta_{1}}$\textsl{. }

\textsl{(1) if }$l(w)<k$\textsl{, then }$\theta_{(\beta_{1},\cdots,\beta_{k}%
)}=0$\textsl{;}

\textsl{(2) if }$l(w)=k$\textsl{, then }$\theta_{(\beta_{1},\cdots,\beta_{k}%
)}$\textsl{ depends only on w and not on the ordered sequence }$\beta
_{1},\cdots,\beta_{k}\in\Delta$\textsl{. In this case we put }$\theta
_{w}=\theta_{(\beta_{1},\cdots,\beta_{k})}$\textsl{.}

\textbf{Proof.} Furnished with properties (1)-(5) in Lemma 4.3 (in analogue
with 3.3 Lemma in [BGG]), an argument parallel to the proof of [BGG, 3.4
Theorem] verifies Lemma 4.4. To do so one needs only to replace $\mathfrak{b}%
_{\mathbb{Q}}^{\ast}=\mathfrak{b}_{\mathbb{Z}}^{\ast}\otimes\mathbb{Q}$ in
[BGG] by $H^{2}(G/T)$, where $\mathfrak{b}_{\mathbb{Z}}^{\ast}$ is the group
of weights of $G$ [Hu. p.67], and to resort to certain properties of Weyl
groups from [BGG; \S2]. For brevity we omit the details (see also Proposition
6 in \textbf{7.5}).

\bigskip

\textbf{4.4.} \textbf{Bott-Samelson cycles and their cohomologies. }Iterating
the construction of leaf spheres gives rise to Bott-Samelson cycles.

\bigskip

For a $y\in G/T$ and an ordered sequence $(\beta_{1},\cdots,\beta_{k})$ of $k$
roots (in which repetitions like $\beta_{i}=\beta_{j}$ for some $1\leq i<j\leq
k$ may occur), we construct an oriented twisted product of $2$-spheres
$S(y;\beta_{1},\cdots,\beta_{k})$ together with a smooth map

$\qquad\quad\varphi_{\beta_{1},\cdots,\beta_{k}}:S(y;\beta_{1},\cdots
,\beta_{k})\rightarrow G/T$

\noindent by induction on $k$. To start with, $S(y;\beta_{1})$ is the leaf
sphere specified by (4.3) with $\varphi_{\beta_{1}}$ the natural inclusion
$S(y;\beta_{1})\subset G/T$. Assume next that the map $\varphi_{\beta
_{1},\cdots,\beta_{k-1}}:S(y;\beta_{1},\cdots,\beta_{k-1})\rightarrow G/T$ has
been defined. Let

\qquad\quad$p_{k-1}:S(y;\beta_{1},\cdots,\beta_{k})\rightarrow S(y;\beta
_{1},\cdots,\beta_{k-1})$

\noindent be the induced bundle of $p_{\beta_{k}}:$ $S(\beta_{k})\rightarrow
G/T$ via $\varphi_{\beta_{1},\cdots,\beta_{k-1}}$. The cross section
$s_{\beta_{k}}$ of $p_{\beta_{k}}$ defines a cross section $s_{k-1}$ of
$p_{k-1}$. The obvious bundle map $\widehat{\varphi}_{\beta_{1},\cdots
,\beta_{k-1}}$ over $\varphi_{\beta_{1},\cdots,\beta_{k-1}}$ followed by
$f_{\beta_{k}}$ gives rise to the desired map

\noindent(4.6) $\qquad\varphi_{\beta_{1},\cdots,\beta_{k}}=f_{\beta_{k}}%
\circ\widehat{\varphi}_{\beta_{1},\cdots,\beta_{k-1}}$.

The final step in the above construction is illustrated by the diagram below.

\noindent(4.7)$\qquad%
\begin{array}
[c]{ccccc}%
S(y;\beta_{1},\cdots,\beta_{k}) & \overset{\widehat{\varphi}_{\beta_{1}%
,\cdots,\beta_{k-1}}}{\rightarrow} & S(\beta_{k}) &  & \\
p_{k-1}\downarrow\uparrow s_{k-1} &  & p_{\beta_{k}}\downarrow\uparrow
s_{\beta_{k}} & \overset{f_{\beta_{k}}}{\searrow} & \\
S(y;\beta_{1},\cdots,\beta_{k-1}) & \overset{\varphi_{\beta_{1},\cdots
,\beta_{k-1}}}{\rightarrow} & G/T & \overset{f_{\beta_{k}}\circ s_{\beta_{k}%
}=id}{\rightarrow} & G/T\text{.}%
\end{array}
$

\bigskip

\textbf{Remark 5.} Alternatively one has

$S(y;\beta_{1},\cdots,\beta_{i})=\{(y_{1},\cdots,y_{i})\in G/T\times
\cdots\times G/T\mid$

$\qquad\qquad\qquad\qquad\qquad y_{1}\in S(y;\beta_{1}),\cdots,y_{i}\in
S(y;\beta_{i})\}$.

\noindent The projection $p_{i-1}$ (resp. the cross section $s_{i-1}$) is
given by

$\qquad p_{i-1}(y_{1},\cdots,y_{i})=(y_{1},\cdots,y_{i-1})$

\qquad(resp. $s_{i-1}(y_{1},\cdots,y_{i-1})=(y_{1},\cdots,y_{i-1},y_{i-1})$).

\noindent The map $\varphi_{\beta_{1},\cdots,\beta_{i}}:S(y;\beta_{1}%
,\cdots,\beta_{i})\rightarrow G/T$ is seen to be

\qquad$\varphi_{\beta_{1},\cdots,\beta_{i}}(y_{1},\cdots,y_{i})=y_{i}$.

\bigskip

\textbf{Definition 6 }(cf. \textbf{7.2}). The map $\varphi_{\beta_{1}%
,\cdots,\beta_{k}}:S(y;\beta_{1},\cdots,\beta_{k})\rightarrow G/T$ is called
\textsl{the Bott-Samelson cycle} associated to the sequence $\beta_{1}%
,\cdots,\beta_{k}$ of roots.

\bigskip

For a subset $L=[i_{1},\cdots,i_{r}]\subseteq\lbrack1,\cdots,k]$, consider the
embedding $i_{L}:S(y;\beta_{i_{1}},\cdots,\beta_{i_{r}})\rightarrow
S(y;\beta_{1},\cdots,\beta_{k})$ by

$\qquad\qquad i_{L}(z_{i_{1}},\cdots,z_{i_{r}})=(y_{1},\cdots,y_{k})$,

\noindent where $y_{t}=z_{i_{s}}$ if $i_{s}\leq t<i_{s+1}$ (cf. Remark 5). Set

$\qquad S(y;L)=\operatorname{Im}\{i_{L}:S(y;\beta_{i_{1}},\cdots,\beta_{i_{r}%
})\rightarrow S(y;\beta_{1},\cdots,\beta_{k})\}$.

\noindent Note that the tower of smooth maps

$\qquad S(y;\beta_{1},\cdots,\beta_{k})\overset{p_{k-1}}{\rightarrow}%
S(y;\beta_{1},\cdots,\beta_{k-1})\overset{p_{k-2}}{\rightarrow}\cdots
\overset{p_{1}}{\rightarrow}S(y;\beta_{1})$,

\noindent together with the sections $s_{i}$, $1\leq i\leq k-1$, furnishes the
space $S(y;\beta_{1},\cdots,\beta_{k})$ with the structure of a twisted
product of $2$-spheres.

\textbf{Lemma 4.5. }\textsl{Let }$[S(y;L)]\in H_{2\mid L\mid}(S(y;\beta
_{1},\cdots,\beta_{k}))$\textsl{ be the fundamental class of the cycle
}$S(y;L)\subset S(y;\beta_{1},\cdots,\beta_{k})$\textsl{. Then}

\textsl{(1) the set }$\{[S(y;L)]\mid L\subseteq\lbrack1,\cdots,k]\}$\textsl{
constitutes an additive basis for the graded }$\mathbb{Z}$\textsl{-module
}$H_{\ast}(M)$\textsl{.}

\noindent\textsl{Further, let} $x_{i}\in H^{2}(S(y;\beta_{1},\cdots,\beta
_{k}))$\textsl{,}$1\leq i\leq k$\textsl{, be the classes Kronecker dual to the
basis }$\{[S(y;i)]\mid$\textsl{ }$1\leq i\leq k\}$ \textsl{of} $H_{2}%
(S(y;\beta_{1},\cdots,\beta_{k}))$\textsl{. Then}

\textsl{(2) the structure matrix of }$S(y;\beta_{1},\cdots,\beta_{k})$
\textsl{(with respect to }$x_{1},\cdots,x_{k}$\textsl{)} \textsl{is
}$A=(a_{ij})_{k\times k}$\textsl{, where }

$\qquad\qquad a_{ij}=\{%
\begin{array}
[c]{c}%
-\beta_{i}\circ\beta_{j}\ \text{\textsl{if} }i<j\text{;}\\
0\ \text{\textsl{if} }i>j\text{\textsl{.}\qquad\quad}%
\end{array}
$

\textsl{(3)} \textsl{the Kronecker pairing in }$S(y;\beta_{1},\cdots,\beta
_{k})$ \textsl{is given by}

$\qquad\qquad<x_{L},[S(y;K)]>=\delta_{L,K}$\textsl{, }$L,K\subseteq
\lbrack1,\cdots,k]$\textsl{,}

\noindent\textsl{where }$x_{L}=x_{i_{1}}\cdots x_{i_{r}}$\textsl{ if }%
$L=$\textsl{ }$[i_{1},\cdots,i_{r}]$\textsl{.}

\textbf{Proof.} (1) and (3) correspond, respectively, to the items (1) and (3)
in Lemma 3.3. It remains to show (2).

For an $2\leq j\leq k$ the normal bundle of the section $s_{j-1}:S(y;\beta
_{1},\cdots,\beta_{j-1})\rightarrow S(y;\beta_{1},\cdots,\beta_{j})$ is seen
to be the the induced bundle $\varphi_{\beta_{1},\cdots,\beta_{j-1}}^{\ast
}E_{\beta_{j}}$, whose Euler class $e_{j}\in H^{2}(S(y;\beta_{1},\cdots
,\beta_{j-1}))$ is $\varphi_{\beta_{1},\cdots,\beta_{j-1}}^{\ast}(e_{\beta
_{j}})$ by the naturality property of Euler classes. That is

$\qquad e_{j}=\varphi_{\beta_{1},\cdots,\beta_{j-1}}^{\ast}(e_{\beta_{j}})$.

\noindent On the other hand, by assuming that the structure matrix of
$S(y;\beta_{1},\cdots,\beta_{k})$ is $A=(a_{ij})_{k\times k}$, we have the expression

$\qquad e_{j}=-a_{1,j}x_{1}-\cdots-a_{j-1,j}x_{j-1}$

\noindent by definition 5, Section 3. Since the $x_{t}$ are Kronecker dual to
the $[S(y;[t])$ we have, for $t<j$, that

$-a_{t,j}=<\varphi_{\beta_{1},\cdots,\beta_{j-1}}^{\ast}(e_{\beta_{j}}),[S(y;[t])]>$

$\qquad=<e_{\beta_{j}},\varphi_{\beta_{1},\cdots,\beta_{j-1}\ast}i_{[t]\ast
}[S(y;\beta_{t})]>$ (by the naturality of $<,>$)

$\qquad=<e_{\beta_{j}},\varphi_{\beta_{t}\ast}[S(y;\beta_{t})]>$ (since
$\varphi_{\beta_{1},\cdots,\beta_{j-1}}\circ i_{[t]}=\varphi_{\beta_{t}}$)

$\qquad=\beta_{t}\circ\beta_{j}$ (by Lemma 4.1).

\noindent This completes the proof of Lemma 4.5.

\section{The induced action of a Bott-Samelson cycle}

Given a sequence $\beta_{1},\cdots,\beta_{k}$ of roots consider the induce
ring map

$\qquad\qquad\varphi_{\beta_{1},\cdots,\beta_{k}}^{\ast}:H^{\ast
}(G/T)\rightarrow H^{\ast}(S(y;\beta_{1},\cdots,\beta_{k}))$

\noindent of the Bott-Samelson cycle $\varphi_{\beta_{1},\cdots,\beta_{k}}%
$\textbf{.} The product in $H^{\ast}(S(y;\beta_{1},\cdots,\beta_{k}))$ is well
understood (cf. Lemma 3.4 and Lemma 4.5). Our aim is to reduce calculations in
the ring $H^{\ast}(G/T)$ (which has been posed to be in question) to that in
$H^{\ast}(S(y;\beta_{1},\cdots,\beta_{k}))$ via $\varphi_{\beta_{1}%
,\cdots,\beta_{k}}^{\ast}$.

Recall from Lemma 4.5 that $H^{\ast}(S(y;\beta_{1},\cdots,\beta_{k}))$ has the
additive basis $\{x_{L}\mid$\textsl{ }$L\subseteq$\textsl{ }$[1,\cdots,k]\}$.
Therefore, for a $u\in H^{2r}(G/T)$, one has a unique expression

\noindent(5.1)$\qquad\qquad\varphi_{\beta_{1},\cdots,\beta_{k}}^{\ast
}(u)=\underset{\mid L\mid=r;L\subseteq\ [1,\cdots,k]}{\Sigma}a_{L}(u)x_{L}$,
$a_{L}(u)\in\mathbb{Z}$.

\noindent The determination of $\varphi_{\beta_{1},\cdots,\beta_{k}}^{\ast}$
amounts to find the $a_{L}(u)$.

\bigskip

In the special cases where $\beta_{1},\cdots,\beta_{k}$ is a sequence of
simple roots, the action of $\varphi_{\beta_{1},\cdots,\beta_{k}}^{\ast}$ on
Schubert classes $P_{w}$, $w\in W$, in $G/T$ can be determined completely.

\textbf{Lemma 5.1.} If $\beta_{1},\cdots,\beta_{k}$ is a sequence of simple
roots, then the induced map $\varphi_{\beta_{1},\cdots,\beta_{k}}^{\ast}$ satisfies

$\qquad\varphi_{\beta_{1},\cdots,\beta_{k}}^{\ast}(P_{w})=(-1)^{l(w)}%
\sum\limits_{\mid L\mid=l(w),\sigma_{L}=w}x_{L}$.

This section is devoted to a proof of Lemma 5.1.

\bigskip

\textbf{5.1.} As the first step we express the coefficient $a_{L}(u)$ in (5.1)
in terms of the divided differences $\theta_{\beta}$. For a subset $L=$
$[i_{1},\cdots,i_{r}]\subseteq\lbrack1,\cdots,k]$ write $\theta_{L}$ for the
composition $\theta_{\beta_{i_{1}}}\circ\cdots\circ\theta_{\beta_{i_{r}}}$.

\textbf{Lemma 5.2. }\textsl{If }$u\in H^{2r}(G/T),$ \textsl{then for
}$L\subseteq\lbrack1,\cdots,k]$\textsl{ with }$\mid L\mid=r$

$\qquad\qquad a_{L}(u)=\theta_{L}(u)$\textsl{ (}$\in H^{0}(G/T)=\mathbb{Z}$\textsl{).}

\textbf{Proof.} In term of the Kronecker pairing one has

$\quad a_{L}(u)=<\varphi_{\beta_{1},\cdots,\beta_{k}}^{\ast}(u),[S(y;L)]>$ (by
(3) of Lemma 4.5)

$=<i_{L}^{\ast}\circ\varphi_{\beta_{1},\cdots,\beta_{k}}^{\ast}(u),[S(y;\beta
_{i_{1}},\cdots,\beta_{i_{r}})]>$ (by the naturality of $<,>$)

$=<\varphi_{\beta_{i_{1}},\cdots,\beta_{i_{r}}}^{\ast}(u),[S(y;\beta_{i_{1}%
},\cdots,\beta_{i_{r}})]>$ (since $\varphi_{\beta_{i_{1}},\cdots,\beta_{i_{r}%
}}=\varphi_{\beta_{1},\cdots,\beta_{k}}\circ i_{L}$).

\noindent This reduces the proof to the special case $L=$\textsl{ }%
$[1,\cdots,k]$. This will be done by induction on $k$. The case $k=1$ is
easily verified by

$\qquad\varphi_{\beta_{1}}^{\ast}(u)=i^{\ast}\circ f_{\beta_{1}}^{\ast
}(u)=i^{\ast}(p_{\beta_{1}}^{\ast}(u)+p_{\beta_{1}}^{\ast}(\theta_{\beta_{1}}(u))x)$

\qquad\qquad$=\theta_{\beta_{1}}(u)x$ (cf. (4.5)).

\noindent Assume, finally, that $L=$\textsl{ }$[1,\cdots,k]$ (i.e. $u\in
H^{2k}(G/T)$). We compute

$\varphi_{\beta_{1},\cdots,\beta_{k}}^{\ast}(u)=\widehat{\varphi}_{\beta
_{1},\cdots,\beta_{k-1}}^{\ast}(f_{\beta_{k}}^{\ast}(u))$ (by (4.6))

\noindent$=\widehat{\varphi}_{\beta_{1},\cdots,\beta_{k-1}}^{\ast}%
(p_{\beta_{k}}^{\ast}(u)+p_{\beta_{k}}^{\ast}(\theta_{\beta_{k}}(u))x)$ (by
the definition of $\theta_{\beta_{k}}$ in (3.3))

\noindent$=p_{k-1}^{\ast}(\varphi_{\beta_{1},\cdots,\beta_{k-1}}^{\ast
}(u))+p_{k-1}^{\ast}(\varphi_{\beta_{1},\cdots,\beta_{k-1}}^{\ast}%
(\theta_{\beta_{k}}(u)))x_{k}$ (by the diagram (4.7))

\noindent$=p_{k-1}^{\ast}(\varphi_{\beta_{1},\cdots,\beta_{k-1}}^{\ast}%
(\theta_{\beta_{k}}(u)))x_{k}$ ( $\varphi_{\beta_{1},\cdots,\beta_{k-1}}%
^{\ast}(u)\in H^{2k}(S(y;\beta_{1},\cdots,\beta_{k-1}))=0$)

\noindent$=p_{k-1}^{\ast}(\theta_{\lbrack1,\cdots,k-1]}(\theta_{\beta_{k}%
}(u))x_{1}\cdots x_{k-1}))x_{k}$ (by the inductive hypothesis)

\noindent$=\theta_{L}(u)x_{L}$.

This finishes the proof.

\bigskip

\textbf{5.2.} \textbf{Bott-Samelson resolution of }$X_{w}$\textbf{.} We refer
to \textbf{7.3} and \textbf{Definition 7} in \textbf{7.4} for two equivalent
geometric descriptions of Schubert varieties $X_{w}$, $w\in W$, in $G/T$. What
is really relevant to us is the desingularization of $X_{w}$, rather than
$X_{w}$ itself, originated from Bott-Samelson in the following way (compare
[BS$_{2}$, p.1000] with discussions in \textbf{7.1}-\textbf{7.4}).

\bigskip

Picture $W$ as the $W$-orbit $\{w(\alpha)\in L(T)\mid w\in W\}$ of the regular
point $\alpha$. For an $w\in W$ let $C_{w}$ be a straight line segment in
$L(T)$ from the Weyl chamber containing $\alpha$ to $w(\alpha)$ that crosses
the planes in $D(G)$ one at a time. Assume that they are met in the order
$L_{\alpha_{1}},\cdots,L_{\alpha_{k}}$, $\alpha_{i}\in\Phi^{+}$. We have
$l(w)=k$ and

$\qquad\qquad w=$ $\sigma_{\alpha_{k}}\circ\cdots\circ\sigma_{\alpha_{1}}$,

\noindent where $\sigma_{\alpha}$ is the reflection of $L(T)$ in $L_{\alpha
}\in D(G)$ (cf. [Han]).

In term of the sequence of positive roots $\alpha_{1},\cdots,\alpha_{k}$
specified by the segment $C_{w}$ we form the composed map

\noindent(5.2)$\qquad\qquad\varphi_{w}=w\circ\varphi_{\alpha_{1},\cdots
,\alpha_{k}}:S(\alpha;\alpha_{1},\cdots,\alpha_{k})\rightarrow G/T$,

\noindent where $\varphi_{\alpha_{1},\cdots,\alpha_{k}}$ is the Bott-Samelson
cycle associated to $\alpha_{1},\cdots,\alpha_{k}$. The map $\varphi_{w}$ may
be appropriately termed as a \textsl{Bott-Samelson resolution of }$X_{w}$ by
the next result shown in [Du$_{2}$, Proposition 3].

\textbf{Lemma 5.3.} \textsl{The map }$\varphi_{w}$\textsl{ is a degree }%
$1$\textsl{ map onto }$X_{w}$\textsl{.}

\bigskip

>From Lemma 5.3 we have, for a Schubert class $P_{w^{\prime}}\in$ $H^{2k}%
(G/T)$, that

$\qquad\varphi_{w}^{\ast}(P_{w^{\prime}})=\delta_{w,w^{\prime}}x_{1}\cdots
x_{k}$.

\noindent On the other hand one has

$\qquad\varphi_{w}^{\ast}(u)=\theta_{\alpha_{1}}\circ\cdots\circ\theta
_{\alpha_{k}}[w^{\ast}(u)]x_{1}\cdots x_{k}$, $u\in H^{2k}(G/T)$

\noindent by Lemma 5.2. These imply that

\textbf{Lemma 5.4.} $\theta_{\alpha_{1}}\circ\cdots\circ\theta_{\alpha_{k}%
}[w^{\ast}(P_{w^{\prime}})]=\delta_{w,w^{\prime}}$.

\bigskip

\textbf{5.3. The operator }$\theta_{w}:H^{2l(w)}(G/T)\rightarrow
H^{0}(G/T)=\mathbb{Z}$\textbf{. }Let $w\in W$ be with $l(w)=k$. Recall from
Lemma 4.4 that if $w=\sigma_{\beta_{1}}\circ\cdots\circ\sigma_{\beta_{k}}$ is
reduced decomposition, then the composition

$\qquad\theta_{w}=\theta_{\beta_{1}}\circ\cdots\circ\theta_{\beta_{k}}%
:H^{\ast}(G/T)\rightarrow H^{\ast-2k}(G/T)$

\noindent is well defined (i.e. not depending on the reduced decomposition
$w=\sigma_{\beta_{1}}\circ\cdots\circ\sigma_{\beta_{k}}$ chosen). This enables
one to evaluate the operator $\theta_{w}$ by using any reduced decomposition
of $w$.

\bigskip

\textbf{Lemma 5.5.} \textsl{With respect to the additive basis }%
$\{P_{w^{\prime}}\mid w^{\prime}\in W,l(w^{\prime})=k\}$\textsl{ of the }%
$2k$\textsl{-dimensional cohomology }$H^{2k}(G/T)$\textsl{, the operator
}$\theta_{w}:H^{2k}(G/T)\rightarrow\mathbb{Z}$\textsl{ is given by }%
$\theta_{w}(P_{w^{\prime}})=(-1)^{l(w)}\delta_{w,w^{\prime}}$\textsl{.}

\textbf{Proof.} Let $\varphi_{w}=w\circ\varphi_{\alpha_{1},\cdots,\alpha_{k}%
}:S(\alpha;\alpha_{1},\cdots,\alpha_{k})\rightarrow G/T$ be the Bott-Samelson
resolution of $X_{w}$ given in (5.2). As in Example 1 we put

\noindent(5.3)$\qquad\beta_{1}=\alpha_{1},\beta_{2}=\sigma_{\alpha_{1}}%
(\alpha_{2}),\cdots,$ $\beta_{k}=\sigma_{\alpha_{1}}\circ\cdots\circ
\sigma_{\alpha_{k-1}}(\alpha_{k})$.

\noindent Then $\beta_{i}\in\Delta$, and $w=\sigma_{\beta_{1}}\circ\cdots
\circ\sigma_{\beta_{k}}$ is a reduced decomposition of $w$. Consequently,

\noindent(5.4)$\qquad\theta_{\beta_{i}}=\sigma_{\alpha_{1}}^{\ast}\circ
\cdots\circ\sigma_{\alpha_{i-1}}^{\ast}\circ\theta_{\alpha_{i}}\circ
\sigma_{\alpha_{i-1}}^{\ast}\circ\cdots\circ\sigma_{\alpha_{1}}^{\ast}$

\noindent by (5.3) and (5) of Lemma 4.3. From (2) of Lemma 4.4 we have

$\qquad\theta_{w}=\theta_{\beta_{1}}\circ\cdots\circ\theta_{\beta_{k}}$

\qquad\qquad$=(\theta_{\alpha_{1}}\circ\sigma_{\alpha_{1}}^{\ast})\circ
(\theta_{\alpha_{2}}\circ\sigma_{\alpha_{2}}^{\ast})\circ\cdots\circ
(\theta_{\alpha_{k}}\circ\sigma_{\alpha_{k}}^{\ast})\circ w^{\ast}$ (by (5.4)).

\noindent Substituting in

$\qquad\theta_{\alpha_{i}}\circ\sigma_{\alpha_{i}}^{\ast}=\theta_{\alpha_{i}%
}\circ(Id-e_{\alpha_{i}}\theta_{\alpha_{i}})$ (by (3) of Lemma 4.3)

$\qquad\qquad=-\theta_{\alpha_{i}}$ (by (1), (4) and (2) of Lemma 4.3)

\noindent we get $\theta_{w}=(-1)^{l(w)}\theta_{\alpha_{1}}\circ\cdots
\circ\theta_{\alpha_{k}}\circ w^{\ast}$. The proof is completed by Lemma 5.4.

\bigskip

\textbf{5.4. Proof of Lemma 5.1. }For an ordered sequence $\{\beta_{1}%
,\cdots,\beta_{k}\}$ of simple roots consider the induced map

$\qquad\varphi_{\beta_{1},\cdots,\beta_{k}}^{\ast}:H^{2r}(G/T)\rightarrow
H^{2r}(S(\alpha;\beta_{1},\cdots,\beta_{k}))$.

\noindent For a $P_{w}\in H^{2l(w)}(G/T)$ we have

$\qquad\varphi_{\beta_{1},\cdots,\beta_{k}}^{\ast}(P_{w})=\underset{\mid
L\mid=l(w);L\subseteq\ [1,\cdots,k]}{\Sigma}\theta_{L}(P_{w})x_{L}$ (by Lemma 5.2)

$\qquad\qquad=\underset{l(\sigma_{L})=\mid L\mid=l(w)}{\Sigma}\theta
_{\sigma_{L}}(P_{w})x_{L}$ (by Lemma 4.4)

$\qquad\qquad=(-1)^{l(w)}\sum\limits_{\mid L\mid=l(w),\sigma_{L}=w}x_{L}$ (by
Lemma 5.5).

\section{Proof of the Theorem}

The Schubert varieties in a generalized flag manifold $G/H$ can be described
by those in $G/T$, where $H$ is the centralizer of a $1$-parameter subgroup in
$G$.

\bigskip

Let $\Delta=\{\beta_{1},\cdots,\beta_{n}\}$ be the set of simple roots
relative to the regular point $\alpha\in L(T)$ (cf. \textbf{2.1}). Assume that
$b\in L(T)\backslash\{0\}$ is a point lying in exactly $d$ of the singular
hyperplanes $L_{\beta_{1}},\cdots,L_{\beta_{n}}$, say $b\in$ $L_{\beta_{1}%
}\cap\cdots\cap L_{\beta_{d}}$. We set $G_{b}$ to be the centralizer of the
$1$-parameter subgroup $\{\exp(tb)\mid t\in\mathbb{R}\}$ in $G$. It is well
known that

(1) \textsl{if }$d=0$\textsl{, then }$G_{b}$\textsl{ is the fixed maximal
torus }$T$\textsl{;}

\noindent and in general

(2) \textsl{every }$H$\textsl{ is conjugated in }$G$\textsl{ to one of the
subgroups }$G_{b}$\textsl{.}

\noindent By (2) we may assume that $H$ is of the form $G_{b}$ for some $b$
taking as the above.

Consequently, $T$ is also a maximal torus of $H$ and the Weyl group
$W^{\prime}$ of $H$, generated by the reflections $\sigma_{\beta_{i}}$,
$k+1\leq i\leq n$, is a subgroup of $W$. As in Section 1, we identify the set
$W/W^{\prime}$ of left cosets of $W^{\prime}$ in $W$ with the subset of $W:$

$\qquad\qquad\overline{W}=\{w\in W\mid I(w^{\prime})\geq I(w)$ for all
$w^{\prime}\in wW^{\prime}\}$.

\bigskip

Consider the standard fibration $p:G/T\rightarrow G/H$. From [BGG, \S5] we have

\textbf{Lemma 6.1.} \textsl{If }$w\in\overline{W}$\textsl{, the map }%
$p$\textsl{ restricts to a degree }$1$\textsl{ map }$X_{w}\rightarrow
X_{w}(H)$\textsl{ between Schubert varieties; and if }$w\notin\overline{W}%
$\textsl{, }$p_{\ast}[X_{w}]=0$\textsl{.}

\textsl{Consequently, the induced map }$p^{\ast}:H^{\ast}(G/H)\rightarrow
H^{\ast}(G/T)$\textsl{ satisfies}

\textsl{\qquad}$p^{\ast}[P_{w}(H)]=P_{w}$\textsl{, }$w\in\overline{W}$\textsl{.}

\bigskip

\textbf{Proof of the Theorem.} For a pair $u,v\in\overline{W}$ assume, in
$H^{\ast}(G/H)$, that

\noindent$\qquad\qquad\qquad P_{u}(H)\cdot P_{v}(H)=\sum\limits_{l(w^{\prime
})=l(u)+l(v),w^{\prime}\in\overline{W}}a_{u,v}^{w^{\prime}}P_{w^{\prime}}(H)$,
$a_{u,v}^{w^{\prime}}\in\mathbb{Z}$.

\noindent Applying the induced ring map $p^{\ast}$ we get the equality

\noindent(6.1)$\qquad$\noindent$\qquad P_{u}\cdot P_{v}=\sum
\limits_{l(w^{\prime})=l(u)+l(v),w^{\prime}\in\overline{W}}a_{u,v}^{w^{\prime
}}P_{w^{\prime}}$, $a_{u,v}^{w^{\prime}}\in\mathbb{Z}$

\noindent in $H^{\ast}(G/T)$ by Lemma 6.1.

Let $w=\sigma_{\beta_{1}}\circ\cdots\circ\sigma_{\beta_{k}}$, $\beta_{i}%
\in\Delta$ be a reduced decomposition of an $w\in\overline{W}$, $k=l(u)+l(v)$,
and let $A_{w}=(a_{i,j})_{k\times k}$ be the associated Cartan matrix of $w$.
Consider the Bott-Samelson cycle $\varphi_{\beta_{1},\cdots,\beta_{k}%
}:S(\alpha;\beta_{1},\cdots,\beta_{k})\rightarrow G/T$ associated to the
sequence $\beta_{1},\cdots,\beta_{k}$ of simple roots. Applying the induced
ring map $\varphi_{\beta_{1},\cdots,\beta_{k}}^{\ast}$ to (6.1) yields in
$H^{\ast}(S(\alpha;\beta_{1},\cdots,\beta_{k}))$ that

$\qquad\varphi_{\beta_{1},\cdots,\beta_{k}}^{\ast}\mathcal{[}P_{u}\cdot
P_{v}]=\sum\limits_{l(w^{\prime})=l(u)+l(v),w^{\prime}\in\overline{W}}%
a_{u,v}^{w^{\prime}}\varphi_{\beta_{1},\cdots,\beta_{k}}^{\ast}[P_{w^{\prime}}]$

\qquad\qquad$=(-1)^{k}a_{u,v}^{w}x_{1}\cdots x_{k}$,

\noindent where the second equality follows from

$\qquad\varphi_{\beta_{1},\cdots,\beta_{k}}^{\ast}[P_{w^{\prime}}]=\{%
\begin{array}
[c]{c}%
(-1)^{k}x_{1}\cdots x_{k}\text{ if }w^{\prime}=w\text{;}\\
0\text{ if }w^{\prime}\neq w\qquad\qquad
\end{array}
$

\noindent by Lemma 5.1. On the other hand

$\qquad\varphi_{\beta_{1},\cdots,\beta_{k}}^{\ast}[P_{u}\cdot P_{v}]$
$=[((-1)^{l(u)}\sum\limits_{\substack{\mid L\mid=l(u)\\\sigma_{L}=u}%
}x_{L})((-1)^{l(v)}\sum\limits_{\substack{\mid K\mid=l(v)\\\sigma_{K}=v}}x_{K})]$

\noindent(again) by Lemma 5.1. Summarizing we get in $H^{2k}(S(\alpha
;\beta_{1},\cdots,\beta_{k}))$ that

$\qquad\qquad(-1)^{k}a_{u,v}^{w}x_{1}\cdots x_{k}=(-1)^{l(u)+l(v)}%
(\sum\limits_{\substack{\mid L\mid=l(u)\\\sigma_{L}=u}}x_{L})(\sum
\limits_{\substack{\mid K\mid=l(v)\\\sigma_{K}=v}}x_{K}).$

\noindent Evaluating both sides on the orientation class $[S(\alpha;\beta
_{1},\cdots,\beta_{k})]$ and noting that $<x_{1}\cdots x_{k},[S(\alpha
;\beta_{1},\cdots,\beta_{k})]>=1$ and $k=l(u)+l(v)$, we get by Lemma 3.4 that

$\qquad\qquad a_{u,v}^{w}=T_{A_{w}}[(\sum\limits_{\substack{\mid
L\mid=l(u)\\\sigma_{L}=u}}x_{L})(\sum\limits_{\substack{\mid K\mid
=l(v)\\\sigma_{K}=v}}x_{K})]$.

\noindent This completes the proof of the Theorem.

\section{Historical remarks}

\textbf{7.1. K-cycles in a flag manifold G/T.} In 1954 R. Bott constructed a
Morse function on $G/T$ and showed that $G/T$ was a cell complex with cells in
the even dimensions only. They turn out to be so-called \textsl{K-cycles} of
Bott-Samelson [BS$_{2}$] formulated in the following plausible way.

For each root $\beta\in\Phi$ let $K_{\beta}\subset G$ be the stabilizer of the
singular plane $L_{\beta}\in D(G)$ under the adjoint action of $G$ on $L(G)$.
For an ordered sequence $\{\beta_{1},\cdots,\beta_{k}\}$ of roots\ one forms
the products $K(\beta_{1},\cdots,\beta_{k})=K_{\beta_{1}}\times\cdots\times
K_{\beta_{k}}$. Since $T\subset K_{\beta_{i}}$ for each $i$, the group
$T(k)=T\times\cdots\times T$ ($k$ factors) acts on $K(\beta_{1},\cdots
,\beta_{k})$ from the right by

$\qquad(c_{1},\cdots,c_{k})\cdot(t_{1},\cdots,t_{k})=(c_{1}t_{1},t_{1}%
^{-1}c_{2}t_{2},\cdots,t_{k-1}^{-1}c_{k}t_{k})$.

\noindent This defines $K(\beta_{1},\cdots,\beta_{k})$ as a $T(k)$-principal
bundle, whose base manifold is called $K_{\beta_{1}}\times_{T}\cdots\times
_{T}K_{\beta_{k}}$. The point in the base corresponding to $(c_{1}%
,\cdots,c_{k})\in K(\beta_{1},\cdots,\beta_{k})$ is denoted by $[c_{1}%
,\cdots,c_{k}]$. The \textsl{K-cycle} associated to the sequence $\{\beta
_{1},\cdots,\beta_{k}\}$ of roots is the map $f_{\beta_{1},\cdots,\beta_{k}%
}:K_{\beta_{1}}\times_{T}\cdots\times_{T}K_{\beta_{k}}\rightarrow G/T$ by

$\qquad\qquad f_{\beta_{1},\cdots,\beta_{k}}[c_{1},\cdots,c_{k}]=Ad_{c_{1}%
\cdots c_{k}}(\alpha)$.

\bigskip

Certain K-cycles were selected to describe the stable manifolds of a perfect
Morse function on $G/T$, hence provide an explicit additive basis for the
homology $H_{\ast}(G/T)$ [BS$_{2}$]. Picture $W$ as the $W$-orbit
$\{w(\alpha)\in L(T)\mid w\in W\}$ of the regular point $\alpha$. For each
$w\in W$ let $C_{w}$ be a straight line segment in $L(T)$ from the Weyl
chamber containing $\alpha$ to $w(\alpha)$ that crosses the planes in $D(G)$
one at a time, and assume that they are met in the order $L_{\alpha_{1}%
},\cdots,L_{\alpha_{k}}$, $\alpha_{i}\in\Phi^{+}$. Let $\Gamma_{w}%
=K_{\alpha_{1}}\times_{T}\cdots\times_{T}K_{\alpha_{k}}$ and define
$g_{w}:\Gamma_{w}\rightarrow G/T$ to be the composition $w\circ f_{\alpha
_{1},\cdots,\alpha_{k}}$. It was shown in [BS$_{2}$] that

\textbf{Proposition 1.} \textsl{The set of cycles }$\{g_{w\ast}[\Gamma_{w}]\in
H_{\ast}(G/T)\mid w\in W\}$\textsl{ is a basis of }$H_{\ast}(G/T)$.

\bigskip

\textbf{7.2. Bott-Samelson cycles in an isoparametric submanifold.} The
embedding $\varphi:G/T\rightarrow L(G)$ given by the adjoint representation at
the beginning of Section 4 defines $G/T$ as an \textsl{isoparametric
submanifolds} in the Euclidean space $L(G)$ [HPT].

In general, associated to any \textsl{isoparametric submanifold} $M$ in an
Euclidean space $\mathbb{R}^{N}$ there are also concepts like (finite) Coxeter
group, root system and Dynkin diagram (marked with multiplicities). In order
to generalize Bott-Samelson's above cited result to such more general spaces
which are also of historical interests in differential geometry,
Hsiang-Palais-Terng introduced in [HPT] the space $S(y;\alpha_{1}%
,\cdots,\alpha_{k})$ as well as the map $\varphi_{\beta_{1},\cdots,\beta_{k}%
}:S(y;\beta_{1},\cdots,\beta_{k})\rightarrow M$ under the name
``\textsl{Bott-Samelson cycles}''in the same way as that given in Remark 5,
Section 4. The construction of these cycles uses only the integrability of the
tangent distributions on $M$ (which are also indexed by sequences of positive
roots relative to a non-focal point $\alpha\in\mathbb{R}^{N}$ of $M$)
\textsl{while the groups }$K_{\beta}$\textsl{'s required to define K-cycles no
long always exist in this more general situation}.

\bigskip

The idea of Bott-Samelson cycles does generalize the K-cycles of Bott-Samelson
in the following sense.

\textbf{Proposition 2 }(cf.[Du$_{2}$, Lemma 8])\textbf{.} \textsl{If }%
$M=G/T$\textsl{ , there is an orientation preserving diffeomorphism
}$g:K_{\beta_{1}}\times_{T}\cdots\times_{T}K_{\beta_{k}}\rightarrow
S(\alpha;\beta_{1},\cdots,\beta_{k})$ \textsl{so that the following mapping
triangle commutes}

$\qquad\qquad\qquad%
\begin{array}
[c]{ccc}%
K_{\beta_{1}}\times_{T}\cdots\times_{T}K_{\beta_{k}} &  & \underset
{f_{\beta_{1},\cdots,\beta_{k}}}{}\\
g\downarrow & \searrow & \\
S(\alpha;\beta_{1},\cdots,\beta_{k}) & \overset{\varphi_{\beta_{1}%
,\cdots,\beta_{k}}}{\rightarrow} & G/T\text{.}%
\end{array}
$

\bigskip

\textbf{7.3. Schubert varieties. }Let $K$ be a linear algebraic group over the
field $\mathbb{C}$ of complex numbers, and let $B\subset K$ be a Borel
subgroup. The homogeneous variety $K/B$ is a projective variety on which the
group $K$ acts by left translations. Historically, Schubert varieties were
introduced in term of the orbits of $B$ action on $K/B$.

Let $T$ be a maximal torus containing in $B$ and let $N(T)$ be the normalizer
of $T$ in $K$. The Weyl group of $K$ (relative to $T$) is $W=N(T)/T$. For an
$w\in W$ take an $n(w)\in N(T)$ such that its residue class mod $T$ is $w$.

The following result was first discovered by Bruhat for classical Lie groups
$K$ in 1954, and proved to be the case for all reductive algebraic linear
groups by Chevalley [Ch$_{2}$] in 1958.

\textbf{Proposition 3.} \textsl{One has the disjoint union decomposition}

$\qquad\qquad K/B=\underset{w\in W}{\cup}Bn(w)\cdot B$

\noindent\textsl{in which each orbit }$Bn(w)\cdot B$\textsl{ is isomorphic to
an affine space of complex dimension }$l(w)$\textsl{.}

The Zariski closure of the open cell $Bn(w)\cdot B$ in $K/B$, denoted by
$X_{w}$, is called the \textsl{Schubert variety associated to }$w$.

\bigskip

\textbf{7.4. Bott-Samelson desingularizations of Schubert varieties.} For a
compact connected Lie group $G$ with a maximal torus $T$ let $K$ be the
complexification of $G$, and let $B$ be a Borel subgroup in $K$ containing
$T$. It is well known that the natural inclusion $G\rightarrow K$ induces an
isomorphism $G/T=K/B$. Conversely, the reductive algebraic linear groups are
exactly the complexifications of the compact real Lie groups (cf. [Ho]).

\bigskip

It follows now from Proposition 1 and 3 that the homology $H_{\ast}(G/T)$ has
two canonical additive bases: the first of these is given by the K-cycles of
Bott-Samelson; the second consists of Schubert varieties, and both of them are
indexed by the Weyl group of $G$.

Before thinking of the problem on the relationship between these two bases of
$H_{\ast}(G/T)$ as ``\textsl{a natural one}'', one should bear in mind that
the first basis was constructed to provide the stable manifolds of a perfect
Morse function on $G/T$, while the second arose from the jumbled efforts
through centuries of great many mathematicians who have contributed to lay the
fundation of algebraic intersection theory for the Schubert's enumerative
calculus [K$_{1}$], [K$_{2}$]. The following result was obtained by Hansen in
1973 [Han].

\textbf{Proposition 4.} \textsl{Under the natural isomorphism }$G/T=K/B$%
\textsl{, the K-cycle }$g_{w}:\Gamma_{w}\rightarrow G/T$\textsl{ of
Bott-Samelson is a degree }$1$\textsl{ map onto the Schubert variety }$X_{w}$\textsl{.}

\bigskip

Combining Proposition 2 with Proposition 4, we have the following alternative
definition of Schubert varieties in $G/T$ (for compact $G$) without resorting
to the complexification of $G$.

Given an $w\in W$ let $C_{w}$ be a straight line segment in $L(T)$ from the
Weyl chamber containing $\alpha$ to $w(\alpha)$ that crosses the planes in
$D(G)$ one at a time. Assume that they are met in the order $L_{\alpha_{1}%
},\cdots,L_{\alpha_{k}}$, $\alpha_{i}\in\Phi^{+}$.

\textbf{Definition 7.} The \textsl{Schubert variety} in $G/T$ associated to an
$w\in W$ is $X_{w}=\operatorname{Im}\varphi_{w}$, where $\varphi_{w}$ is the composition

$\qquad\qquad\varphi_{w}=w\circ\varphi_{\alpha_{1},\cdots,\alpha_{k}}%
:S(\alpha;\alpha_{1},\cdots,\alpha_{k})\rightarrow G/T$.

This is the version of the descriptions of Schubert varieties that we have
made use of in Section 5 and 6.

\bigskip

\textbf{7.5. Divided differences.} The inclusion $T\rightarrow G$ of the
maximal torus induces a fibration

$\qquad\qquad G/T\overset{c}{\hookrightarrow}BT\rightarrow BG$,

\noindent where $BH$ denotes the classifying space of a Lie group $H$ , and
where the fibre inclusion $c$ is equivariant with respect to the standard
$W$-action on both $G/T$ and $BT$ [B].

Let $H^{\ast}(X;\mathbb{R})$ be the cohomology ring (resp. algebra) of a space
$X$ with coefficients in the field $\mathbb{R}$ of reals. We recall the
classical result due to Borel [B].

\textbf{Proposition 5.} \textsl{The map }$c$\textsl{\ induces an }${W}%
$\textsl{-equivariant surjective homomorphism of algebras}

\textsl{\ \qquad}$c^{\ast}:{H}^{\ast}{(BT;}\mathbb{R}{)\rightarrow H}^{\ast
}{(G/T;\mathbb{R})}$

\noindent\textsl{with kernel }${H}^{+}{(BT;\mathbb{R})}^{W}$\textsl{, the
ideal in }${H}^{\ast}{(BT;\mathbb{R})}$ \textsl{generated by }$W$%
\textsl{-invariants in positive degrees.}

The induced map $c^{\ast}$, playing a key role in this result, is well known
as the \textsl{Borel's characteristic map}.

\bigskip

The infinite complex projective space $\mathbb{C}P^{\infty}$ serves both as
the classifying space of the circle group $S^{1}$ and the Eilenberg-MacLane
space $K(\mathbb{Z},2)$. Keeping this in mind one get two ingredients from
each root $\alpha\in\Phi$.

(1) The reflection $\sigma_{\alpha}$ on $L(T)$ in the hyperplane $L_{\alpha}$
preserves the \textsl{unit lattice }$\Lambda=\exp^{-1}(e)$, hence induces in
successive manner an automorphism of the torus $T=L(T)/\Lambda$; a
diffeomorphism of the classifying space $BT$, and finally, an induced
automorphism $\sigma_{\alpha}^{\ast}$ of the ring $H^{\ast}(BT;\mathbb{R})$.

(2) the \textsl{co-root} $\alpha^{\ast}:L(T)\rightarrow\mathbb{R}$ related to
the root $\alpha\in\Phi$ in the fashion $\alpha^{\ast}(v)=(\alpha_{i},v)$,
$v\in L(T)$ satisfies $\alpha^{\ast}(\Lambda)\subset\mathbb{Z}$, hence induces
successively a homomorphism $T\rightarrow S^{1}$, a map between classifying
space $BT\rightarrow BS^{1}=\mathbb{C}P^{\infty}$, an finally a $2$-cocycle
$[\alpha]\in H^{2}(BT;\mathbb{R})$.

\noindent In terms of these, an additive operation $A_{\alpha}:$ $H^{\ast
}(BT;\mathbb{R})\rightarrow H^{\ast}(BT;\mathbb{R})$ of degree $-2$ can be
defined as follows.

$\qquad\qquad A_{\alpha}(f)=\frac{f-\sigma_{\alpha}^{\ast}(f)}{[\alpha]}$,
$f\in H^{\ast}(BT;\mathbb{R})$.

\bigskip

The operator $A_{\alpha}$ on $H^{\ast}(BT;\mathbb{R})$, known as \textsl{the
divided difference operators associated to the root }$\alpha$, was introduced
independently by Bernstein et al [BGG] and Demazure [De$_{1}$] in 1973. These
operators possess the same properties as that of the operators $\theta
_{\alpha}$ on ${H}^{\ast}{(G/T)}$ (compare [BGG, 3.3 Lemma] with (2) of Lemma
4.3 in Section 4). In fact, Borel's characteristic map gives rise to a linkage
between these two sets of operators.

\textbf{Proposition 6.} \textsl{For each }$\alpha\in\Phi$\textsl{, the two
operators }$A_{\alpha}$\textsl{ and }$\theta_{\alpha}$\textsl{ satisfy the
commutative diagram}

$\qquad%
\begin{array}
[c]{ccc}%
H^{\ast}(BT;\mathbb{R}) & \overset{c^{\ast}}{\rightarrow} & H^{\ast
}(G/T,\mathbb{R})=H^{\ast}(G/T)\otimes\mathbb{R}\\
A_{\alpha}\downarrow &  & \downarrow\theta_{\alpha}\otimes1\\
H^{\ast}(BT;\mathbb{R}) & \overset{c^{\ast}}{\rightarrow} & H^{\ast
}(G/T,\mathbb{R})=H^{\ast}(G/T)\otimes\mathbb{R}\text{.}%
\end{array}
$

\bigskip

However, by introducing the operators $\theta_{\alpha}$ in the much more
general context of spherical represented involutions (cf. \textbf{3.2}), and
by developing their basic properties (Lemma 3.2) in a way independent of the
$A_{\alpha}$ and the Borel characteristic map, these operators not only can
act directly on integral ring ${H}^{\ast}{(G/T)}$ (in which we are actually
interested), but also become applicable in the study of intersection theory of
isoparametric submanifolds, transnormal submanifolds, as well as their focal
manifolds [HPT], [R]. In the context of Morse theory, these families of
manifolds are seen as natural generalizations of flag manifolds in the sense
that \textsl{generalized Schubert cycles} (cf. [HPT, p.449]) can be
constructed for distance functions on these manifolds and which are also
indexed by the cosets of certain finite groups. In these general settings
there are no analogue of Cartan subalgebra (and consequently, Borel's
characteristic maps) available.

\bigskip

\textbf{7.6. Positively multiplying Schubert classes} (continuing from Remark
1, Section 2). Let $S$ be a twisted products of $2$-spheres of rank $k$ with
structure matrix $A$, and let $B=\{x_{I}\mid I\subset\lbrack1,\cdots,k]\}$ be
the additive basis for $H^{\ast}(S)$ specified by Lemma 3.3. Product among the
basis elements of complementary dimensions yields the numbers $c_{I,J}%
\in\mathbb{Z}$ with

$\qquad\qquad x_{I}\cdot x_{J}=c_{I,J}x_{1}\cdots x_{k}$, $\mid I\mid+\mid
J\mid=k$.

\noindent Equivalently, $c_{I,J}=T_{A}(x_{I}\cdot x_{J})$.

The advantages to work with the basis $B$ are obvious. It consists of the
elements Kronecker dual to the ready-made geometric cycles $S(I)$ in $S$
(Lemma 3.3). With respect to this basis the action of Bott-Samelson cycles
admits a simple description (Lemma 5.1). As a result the geometric essence of
the Theorem is transparent. On the other hand, the product in $S$ between
these basis elements are not always positive in the sense that $c_{I,J}<0$ may
occur (cf. Example 3). This explains \ the reason that our formula falls short
of positivity.

\bigskip

The classical Richardson-Littlewood rule for multiplying Schubert classes in
the Grassmannian has the merit to meet the standard of positivity. However, it
is difficult to summarize the rule into an explicit formula required by
effective computation. Its lengthy and technical statement (cf. [St, p.228])
form a sharp contrast with the conciseness of the Theorem.

Concerning the geometric origin of the problem, effective computability should
be granted with the first priority among various standards on the
multiplicative rules [K$_{1}$], [K$_{2}$]. Therefore, it is natural to make
inquires about a unified solution to the problem in its deserved simplicity
and natural generality.

\bigskip

\begin{center}
\textbf{References}
\end{center}

[B] A. Borel, Topics in the homology theory of fiber bundles, Berlin,
Springer, 1967.

[BGG] I. N. Bernstein, I. M. Gel'fand and S. I. Gel'fand, Schubert cells and
cohomology of the spaces G/P, Russian Math. Surveys 28 (1973), 1-26.

[Bi] S. Billey, Kostant polynomials and the cohomology ring for G/B, Duke J.
Math. 96, No.1(1999), 205-224.

[BH] S. Billey and M. Haiman, Schubert polynomials for the classical groups,
Journal of AMS, 8(2)(1995), 443-483.

[BS$_{1}$] R. Bott and H. Samelson, The cohomology ring of G/T, Nat. Acad.
Sci. 41 (7) (1955), 490-492.

[BS$_{2}$] R. Bott and H. Samelson, Application of the theory of Morse to
symmetric spaces, Amer. J. Math., Vol. LXXX, no. 4 (1958), 964-1029.

[CF] I. Ciocan-Fontanine, On quantum cohomology rings of partial flag
varieties, Duck Math. J., Vol.98, No. 3(1999), 485-524.

[Ch$_{1}$] C. Chevalley, La th\'{e}orie des groupes alg\'{e}briques, Proc.
1958 ICM, Cambridge Univ. Press, 1960, 53-68.

[Ch$_{2}$] C. Chevalley, Sur les D\'{e}compositions Celluaires des Espaces
G/B, in Algebraic groups and their generalizations: Classical methods, W.
Haboush ed. Proc. Symp. in Pure Math. 56 (part 1) (1994), 1-26.

[De$_{1}$] M. Demazure, Invariants sym\'{e}triques entiers des groupes de Weyl
et torsion, Invent. Math. 21 (1973), 287-301.

[De$_{2}$] M. Demazure, D\'{e}singularization des vari\'{e}t\'{e}s de Schubert
g\'{e}n\'{e}ralis\'{e}es, Ann. Sci. \'{E}cole. Norm. Sup. (4) 7(1974), 53-88.

[Du$_{1}$] H. Duan, Some enumerative formulas on flag varieties, Communication
in Algebra, 29 (10) (2001), 4395-4419.

[Du$_{2}$] H. Duan, The degree of a Schubert variety, to appear in Adv. in
Math (Available online: http://www.sciencedirect.com /science/journal/00018708).

[DZ$_{1}$] H. Duan and Xuezhi Zhao, Steenrod operations on Schubert classes,
arXiv: math.AT/0306250.

[DZ$_{2}$] H. Duan and Xuezhi Zhao, A program for multiplying Schubert
classes, in preparation.

[DZZ] H. Duan, Xu-an Zhao and Xuezhi Zhao, Cartan matrix and enumerative
calculus, preprint.

[FP] W. Fulton and P. Pragacz, Schubert Varieties and Degeneracy Loci, LNM
1689, Springer, 1998.

[Han] H.C. Hansen, On cycles in flag manifolds, Math. Scand. 33 (1973), 269-274.

[Ho] G. Hochschild, The structure of Lie groups, Holden-Day, Inc., San Francis

\noindent co-London-Amsterdam, 1965.

[HPT] W. Y. Hsiang, R. Palais and C. L. Terng, The topology of isoparametric
submanifolds, J. Diff. Geom., Vol. 27 (1988), 423-460.

[Hu] J. E. Humphreys, Introduction to Lie algebras and representation theory,
Graduated Texts in Math. 9, Springer-Verlag New York, 1972.

[K$_{1}$] S. Kleiman, Problem 15. Rigorous fundation of the Schubert's
enumerative calculus, Proceedings of Symposia in Pure Math., 28 (1976), 445-482.

[K$_{2}$] S. Kleiman, Intersection theory and enumerative geometry: A decade
in review, Algebraic Geometry, Bowdoin 1985 (Spancer Bloch, ed.), Proc.
Sympos. Pure Math., vol. 46, Part 2, AMS. 1987, 321-370.

[KK] B. Kostant and S. Kumar, The nil Hecke ring and the cohomology of G/P for
a Kac-Moody group G, Adv. Math. 62(1986), 187-237.

[LR] D. E. Littlewood and A. R. Richardson, Group characters and algebra,
Philos. Trans. Roy. Soc. London. 233(1934), 99-141.

[LS] A. Lascoux and M.-P. Sch\"{u}tzenberger, Polyn\^{o}mes de Schubert, C. R.
Acad. Sci. Paris 294(1982), 447-450.

[M] I. G. Macdonald, Symmetric functions and Hall polynomials, Oxford
Mathematical Monographs, Oxford University Press, Oxford, second ed., 1995.

[R] S. A. Robertson, Smooth curves of constant width and transnormality, Bull.
London Math. Soc. 16, no.3(1984), 264-274.

[S] F. Sottile, Four entries for Kluwer encyclopaedia of Mathematics, arXiv:
Math. AG/0102047.

[St] R. Stanley, Some combinatorial aspects of the Schubert calculus,
Combinatoire et repr\'{e}sentation du groupe sym\'{e}trique, Strasbourg
(1976), 217-251.
\end{document}